\documentclass[reqno]{article}

\usepackage{orcidlink}
\usepackage{amsmath, amssymb, amsthm}
\usepackage{mathtools}
\usepackage{graphicx}
\usepackage{hyperref}
\usepackage{tikz}
\usepackage{tikz-cd}
\usepackage{ytableau}
\usepackage{enumitem}
\usepackage{algorithm}
\usepackage{algorithmic}
\usepackage{placeins}
\usepackage{multirow}
\usepackage{mathrsfs}
\usetikzlibrary{patterns, shadows}
\usepackage{cite}

\usepackage{bbding}            
\usepackage{xcolor} 
\definecolor{yellow}{HTML}{F0E442} 
\definecolor{green}{HTML}{56B4E9} 
\definecolor{orange}{HTML}{D55E00} 

\usepackage{bm}

\usepackage{newunicodechar}
\usepackage{ytableau}
\usepackage{ytableau,xcolor}
\usepackage{tcolorbox}
\usepackage{mathtools} 
\usepackage{bm}  
\usepackage{stmaryrd}
\usepackage{url}
\usepackage{colortbl} 
\usepackage{authblk}   

\usepackage[
  left=2.5cm,
  right=2.5cm,
  top=2.5cm,
  bottom=2.5cm
]{geometry}

\tcbuselibrary{skins, breakable}
\usetikzlibrary{positioning}
\usetikzlibrary{matrix, arrows.meta} 

\makeatletter
\newcommand{\fixed@sra}{$\vrule height 2\fontdimen22\textfont2 width 0pt\shortrightarrow$}
\newcommand{\shortarrow}[1]{%
  \mathrel{\text{\rotatebox[origin=c]{\numexpr#1*45}{\fixed@sra}}}
}
\makeatother

\makeatletter
\newcommand{\ostar}
{\mathbin{\mathpalette\make@circled\star}}
\newcommand{\make@circled}[2]{%
  \ooalign{$\m@th#1\smallbigcirc{#1}$\cr\hidewidth$\m@th#1#2$\hidewidth\cr}%
}
\newcommand{\smallbigcirc}[1]{%
  \vcenter{\hbox{\scalebox{0.77778}{$\m@th#1\bigcirc$}}}%
}
\makeatother

\DeclareMathOperator{\lcm}{lcm}
\DeclareMathOperator{\conj}{conj}
\DeclareMathOperator{\rev}{rev}
\DeclareMathOperator{\diag}{diag}

\newunicodechar{−}{\ensuremath{-}}

\newtheorem{theorem}{Theorem}[section]
\newtheorem{definition}[theorem]{Definition}
\newtheorem{lemma}[theorem]{Lemma}
\newtheorem{corollary}[theorem]{Corollary}
\newtheorem{proposition}[theorem]{Proposition}

\newtheorem{remark}[theorem]{Remark}
\newtheorem{example}[theorem]{Example}

\title{Bulgarian solitaire with distinguishable cells}

\author{Vikum Bandara\thanks{Email: vbandar@clemson.edu} }
\author{Svetlana Poznanovi\'c\thanks{Email: spoznan@clemson.edu}}

\affil{School of Mathematical and Statistical Sciences\\
Clemson University\\
Clemson, SC 29634, USA}

\date{}

\begin{document}

\maketitle

\begin{abstract}
Classical Bulgarian solitaire tracks only the sizes of the stacks. We ask how the lengths of its periodic orbits change when individual cards are distinguished. To study this question, we introduce four actions $T_1,T_2,T_3,T_4$ on partitions with labeled cells. The actions differ in whether the top or bottom cell is removed from each column and in the order in which the removed cells are arranged to form a new column. We determine the exact sizes of the periodic orbits for all four actions. The main work lies in the analysis of $T_2$, for which we develop a coloring algorithm on Ferrers diagrams. The coloring encodes the cycle structure of the permutation of the labels induced by one full period of the underlying partition.
\end{abstract}

\section{Introduction}
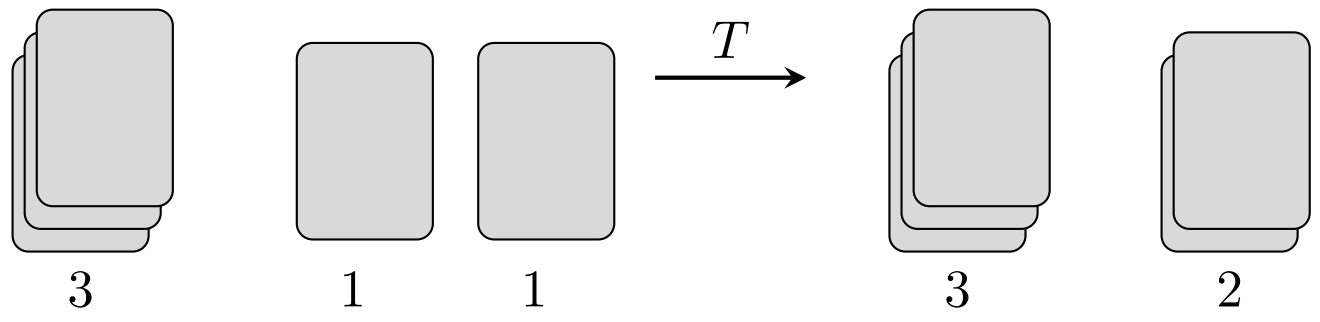
\begin{figure}[h]
\centering
\begin{tikzpicture}[
    card/.style={
        rectangle,
        draw=black,
        rounded corners=3pt,
        minimum width=0.9cm,
        minimum height=1.3cm,
        fill=gray!30,          
        font=\small\bfseries,
        align=center
    }
]

\node[card] at (0,0) {};                
\node[card, xshift=0.08cm, yshift=0.15cm] at (0,0) {};
\node[card, xshift=0.16cm, yshift=0.30cm] at (0,0) {};   

\node[card, xshift=0.08cm, yshift=0.08cm] at (1.8,0) {};

\node[card, xshift=0.08cm, yshift=0.08cm] at (3.0,0) {};

\node at (0,-0.9) {3};
\node at (1.8,-0.9) {1};
\node at (3.0,-0.9) {1};

\draw[->, thick, >=stealth] (3.8,0.5) -- (4.8,0.5) node[midway, above] {$T$};

\begin{scope}[shift={(5.8,0)}]
    \node[card] at (0,0) {};
    \node[card, xshift=0.08cm, yshift=0.15cm] at (0,0) {};
    \node[card, xshift=0.16cm, yshift=0.30cm] at (0,0) {};
    \node at (0,-0.9) {3};

    \node[card] at (1.8,0) {};
    \node[card, xshift=0.08cm, yshift=0.15cm] at (1.8,0) {};
    \node at (1.8,-0.9) {2};
\end{scope}

\end{tikzpicture}
\caption{One move of Bulgarian solitaire from $(3,1,1)$ to $(3,2)$. Remove one card from each stack to create a stack of size 3; the original stack of size 3 becomes a stack of size 2, and the two singleton stacks disappear.}
\label{fig:bulgarian-facedown}
\end{figure}

 Bulgarian solitaire is a mathematical game played with cards. One starts with $n$ cards divided into stacks and placed face down. One step of the game consists of taking one card from every stack to form a new stack.  Figure~\ref{fig:bulgarian-facedown} shows one move in a game with 5 cards, from the configuration with stack sizes $(3,1,1)$ to the configuration $(3,2)$. Since the cards are treated as indistinguishable, a configuration is encoded by its stack sizes. These stack sizes are traditionally ordered in weakly decreasing order and represented as a partition of $n$. Thus, Bulgarian solitaire defines a discrete dynamical system
on the finite set $\mathcal{P}_n$ of partitions of $n$, and every
sequence of iterates is eventually periodic.

According to the survey~\cite{hopkins2012}, which also discusses
the origin of the game's name, Bulgarian solitaire initially
spread among mathematicians through informal conversations.
A central observation was that, when the number of cards is
triangular, $n=k(k+1)/2$, every initial configuration eventually
settles into the partition $(k,k-1,\dots,1)$. This observation
inspired several mathematicians to seek a proof (see, e.g.,~\cite{bojanov1981problem,eriksson1981bulgarisk}) and to investigate the behavior for general $n$. Early analyses of the general case appeared in~\cite{toom1981,brandt1982cycles};
see also~\cite{akinDavis1985bulgarian,hart2012}.
Gardner popularized the game in his \emph{Scientific American}
column, bringing it to an international
audience~\cite{gardner1997bulgarian}.


The dynamics of Bulgarian solitaire have been studied from several
perspectives, including the enumeration of periodic orbits and
the study of periodicity~\cite{brandt1982cycles,igusa1985,griggs1998,hart2012},
as well as the analysis of Garden of Eden
partitions~\cite{hopkins2006,hopkins2007garden}.
Related work examines the game graph and the distribution of
level sizes~\cite{erikssonJonsson2017levels,HarrisNguyen2025Depth,pham2023limiting}.
Beyond the classical game, numerous variants have been studied
(see, e.g.,~\cite{MR4011610,canningsHaigh1992montreal,
olson2016variants,erikssonJonssonSjostrand2020limit,hopkins2024henrik,
yeh1995remarkable,hopkins2009envelopes,meyerowitzEtAl2026floridian}),
including probabilistic versions
(see, e.g.,~\cite{popov2005,erikssonJonssonSjostrand2018randomq}).

In the classical game, the cards are treated as indistinguishable
and only the stack sizes are tracked. This paper asks what happens
to individual cards as the game proceeds. A return of the stack
sizes need not return each card to its original position, so
the labeled configurations can have longer periods than their
underlying partitions. To study this behavior, we extend Bulgarian
solitaire to partitions with distinguishable cells, denoted by
$\mathcal{P}_n^\star$, and define four natural actions
$T_1,T_2,T_3,T_4$. These actions differ in whether a card is
removed from the top or the bottom of each stack and in the
order in which the removed cards are arranged to form a new stack. We determine the sizes of the periodic orbits for all four actions.
The main work lies in the analysis of $T_2$, for which we develop
a coloring algorithm on Ferrers diagrams that describes how the
cards are permuted over one period of the underlying partition.

The paper is organized as follows. Section~\ref{sec:background} reviews  partitions, Ferrers diagrams, and classical Bulgarian solitaire, including the characterization of periodic points. Section~\ref{sec:actions} introduces the four actions, establishes
their relationship with the classical game, and uses reversal
to reduce the orbit-size analysis to $T_1$ and $T_2$. Section~\ref{sec:orbit sizes} determines the periodic orbit sizes
under $T_1$ and develops the coloring algorithm for $T_2$. Section~\ref{subsec:orbit-size} completes the analysis by
determining the exact periodic orbit sizes under $T_2$.
Finally, Section~\ref{sec:conjugation} studies the interaction
with conjugation and the resulting orbit-reversing symmetries.

\section{Background}
\label{sec:background}

In this section, we recall the necessary definitions of partitions, Ferrers diagrams, and classical Bulgarian solitaire, in particular the characterization of periodic points, which will be used to analyze the generalization introduced later.

\subsection{Partitions and Ferrers diagrams}
\begin{definition}
A partition of a positive integer \(n\) is a finite sequence \(\pi = (\pi_1, \pi_2, \dots, \pi_k)\) of positive integers such that $\pi_1 \ge \pi_2 \ge \cdots \ge \pi_k  \quad \text{and} \quad \sum_{i=1}^k \pi_i = n.$ The integers \(\pi_i\) are called the parts of the partition, and \(k = \ell(\pi)\) is called the length of the partition.  The set of all partitions of \(n\) is denoted by \(\mathcal{P}_n\).
\end{definition}

\begin{figure}[H]
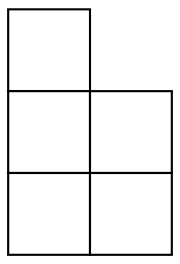


\[ \begin{ytableau} 
*(white)  \\ 
*(white) & *(white) \\ 
*(white) & *(white) 
\end{ytableau} \]

\caption{Rotated Ferrers diagram of $\pi=(3,2)\in\mathcal{P}_{5}$.}
\label{pi_in_P}
\end{figure}

A partition \(\pi = (\pi_1, \pi_2, \dots ,\pi_k)\) of \(n\) is traditionally represented by its Ferrers diagram: a left-aligned array of cells in which the \(i\)-th row from the top contains \(\pi_i\) cells. In our pictorial representation, we rotate this diagram by \(90^\circ\) counterclockwise so that the parts correspond to columns (representing stacks of cards) rather than rows (e.g., Figure~\ref{pi_in_P}). We refer to the cells of this diagram using Euclidean coordinates, with the bottom-left cell at position \((1,1)\). 

For two partitions $\mu$ and $\pi$, we define their union $\mu \cup \pi$ to be the partition obtained by combining the parts of $\mu$ with those of $\pi$ and arranging the resulting entries in weakly decreasing order. For example, $(3,2) \cup (5,4, 2, 1) = (5, 4, 3, 2, 2, 1).$

\subsection{Classical Bulgarian solitaire}

Bulgarian solitaire is a discrete dynamical system defined on $\mathcal{P}_n$. The system evolves through an iterative process governed by the following rule: from a given partition, remove one cell from each part (column) and use the removed cells to form a new part. Formally, Bulgarian solitaire is the map $T: \mathcal{P}_n \to \mathcal{P}_n$ defined by
\[
T(\pi) = (\pi_1 - 1, \pi_2 - 1, \dots, \pi_{\ell(\pi)} - 1) \cup (\ell(\pi)),
\]
where parts of size zero are omitted.

The set
\[
O_T(\pi)=\{\pi,T(\pi),T^2(\pi),\ldots\}
\]
is called the orbit of $\pi$ under $T$. Since $\mathcal{P}_n$ is finite,
the sequence of iterates is eventually periodic.

A partition \( \pi \) is said to be \emph{periodic} under the transformation \( T \) if there exists a positive integer \( \ell \) such that
\[
T^\ell (\pi) = \pi.
\]
The smallest such integer \( \ell \) is called the \emph{period} of \( \pi \) under \( T \).

\begin{definition}
For a fixed positive integer $k$, the $k$-th upper diagonal, denoted by $\diag[k]$, is the set of positions
\[
\diag[k] = \{(1,k), (2,k-1), \ldots, (k,1)\}.
\]
\end{definition}

For a positive integer $n$, let $m$ be the unique integer such that the $m$-th triangular number $\Delta_m = \frac{m(m+1)}{2}$ satisfies
\[
\Delta_m \le n < \Delta_{m+1}.
\]
Define $m' = n - \Delta_m$.

The periodic points of $T$ were characterized by Toom as follows.

\begin{theorem}[\cite{toom1981, hart2012}]
\label{periodic points under T}

Let \(n\ge 1\), let \(m\) be the unique integer such that $\Delta_m \le n < \Delta_{m+1}$ and let $m' = n - \Delta_m$. Let \(\pi\) be a partition of \(n\), and let \(F\) denote its Ferrers diagram. Then \(\pi\) is a periodic point of \(T\) if and only if
\[
|F \cap \diag[k]| =
\begin{cases}
k, & 1 \le k \le m, \\[4pt]
m', & k = m+1, \\[4pt]
0, & k > m+1.
\end{cases}
\]
\end{theorem}

In the special case where $n = \Delta_m$, $\mathcal{P}_n$ contains only one periodic point, namely,
\[
(m, m - 1, m - 2, \dots, 2, 1).
\]

When $n$ is not triangular, the period $\ell$ of a periodic partition divides $m+1$. As can be gleaned from Example~\ref{ex:OT(pi311)}, for periodic points, the newly formed part is always at least as large as all remaining parts and can be placed in the first position. The overall change in the Ferrers shape amounts to a cyclic downward shift of the cells in $\diag[m+1]$. 

The number of distinct periodic orbits of $\mathcal{P}_n$ under the action of $T$ is equal to the number of necklaces consisting of $m^\prime$ black beads and $m + 1 - m^\prime$ white beads, where the symmetry group is the cyclic group of order $m + 1$. Their number is given by the following elegant formula due to Brandt.

\begin{theorem}[\cite{brandt1982cycles}]\label{theoremdistinctcycles}
Given $n \geq 1$, write $n$ uniquely as $\Delta_m + m'$ for some integer $m'$ with $0 \le m' \le m$. The number of distinct periodic orbits in $\mathcal{P}_n$ under $T$ is
\[
\frac{1}{m + 1} \sum_{d \mid \gcd(m + 1, m')} \phi(d) \binom{\frac{m + 1}{d}}{\frac{m'}{d}},
\]
where $\phi$ is Euler’s totient function.
\end{theorem}

\begin{example}\label{ex:OT(pi311)}

Figure~\ref{OT(pi311)} illustrates the orbit of $\pi = (3,1,1)$ under $T$. By Theorem~\ref{periodic points under T}, $\pi$ is a periodic point of $T$. The action of $T$ corresponds to cyclically shifting the cells in $\diag[m+1]$ down the diagonal.

\begin{figure}[H]
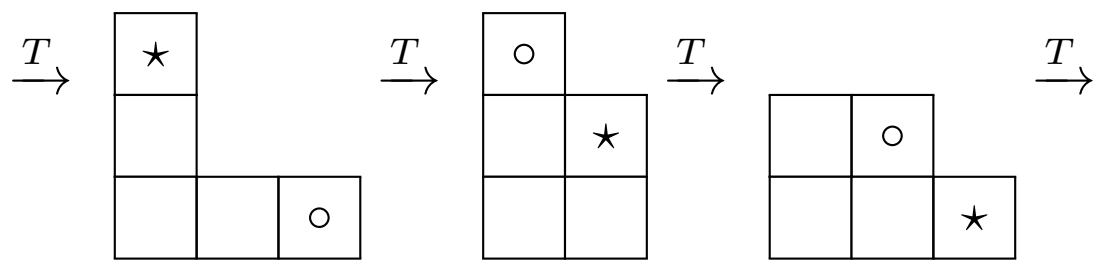

\[
\xrightarrow{T}\,\,\, 
\begin{ytableau} \star  \\ *(white) \\ *(white) & *(white) & \circ \end{ytableau} 
\xrightarrow{T}\,\,\, 
\begin{ytableau} \circ \\ *(white) & \star \\ *(white) & *(white)  \end{ytableau}
\xrightarrow{T}\,\,\, 
\begin{ytableau} \none & \none & \none \\*(white) & \circ \\ *(white) & *(white) & \star \end{ytableau}
\xrightarrow{T}\,\,\,   
\]
\caption{Orbit of $\pi = (3,1,1)$ under $T$.}
\label{OT(pi311)}
\end{figure}

In this case, $n=5=\Delta_2+2$, so $m=2$ and $m'=2$. The period is $\ell=3$, which divides $m+1$. This orbit corresponds to the unique necklace with $m'=2$ black beads and
$m+1-m'=1$ white bead, shown in Figure~\ref{necklace}.

\begin{figure}[H]
\centering
\begin{tikzpicture}
    \draw[gray, thick] (0,0) circle (1.5cm);
    \foreach \angle/\color in {90/white, -30/black, -150/black} {
        \node[circle, fill=\color, draw=black, minimum size=6mm, line width=0.5pt] at (\angle:1.5cm) {};
    }
\end{tikzpicture}
\caption{A necklace with $2$ black beads and $1$ white bead.}
\label{necklace}
\end{figure}
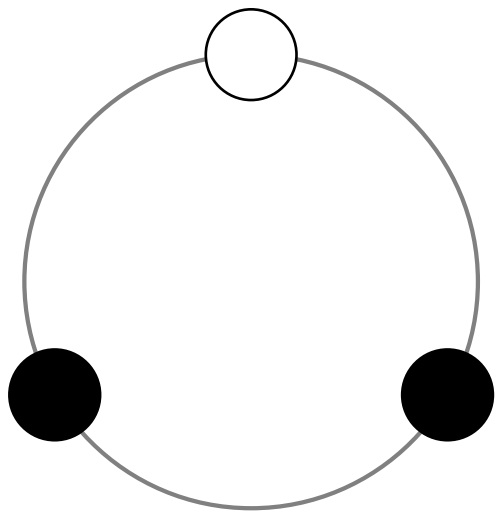

\end{example}

The \emph{depth} $D(n)$ of classical Bulgarian solitaire is the
maximum, over $\pi\in\mathcal{P}_n$, of the minimum number of
applications of $T$ needed to reach a periodic partition. Igusa~\cite{igusa1985} and Etienne~\cite{etienne1991} proved that
$D(\Delta_m)=m(m-1)$ for $m\geq1$. Griggs and Ho~\cite{griggs1998}
gave an upper bound for non-triangular numbers:
$D(\Delta_m+m')\leq m^2-2$ for $m\geq3$ and $1\leq m'\leq m$,
with equality when $m'=m$. For results on counting partitions according to the number of moves required to reach a given periodic orbit, see~\cite{HarrisNguyen2025Depth,pham2023limiting}.

\section{\texorpdfstring
  {Four variants of Bulgarian solitaire on $\mathcal{P}_n^\star$}
  {Four variants of Bulgarian solitaire on partitions with distinguishable cells}} \label{sec:actions}

When the cards are turned face up, each card becomes individually identifiable, so each stack is represented by an ordered sequence of distinct labels. This motivates the study of \emph{partitions with distinguishable cells}, where the identity of each card is taken into account. 

A partition with distinguishable cells is a Ferrers diagram whose cells are bijectively labeled by a fixed set of $n$ labels. No order conditions are imposed on the labels. We write each part as $\pi_i=(\pi_{i1},\pi_{i2},\ldots,\pi_{in_i})$,
where $n_i$ is the length of the $i$-th column and the entries represent distinct cells. We depict each $\pi_i$ as a column in the Ferrers diagram, and $\pi_{ij}$ denotes the labeled cell in position $(i,j)$, that is, the cell in the $i$-th column and $j$-th row. 

We denote the set of all partitions of $n$ with distinguishable cells by $\mathcal{P}_n^\star $. The \textit{set of all labeled cells} of a partition's Ferrers diagram is denoted by $\mathcal{L}_\pi$.

\begin{example}\label{pistar}
Consider the partition with distinguishable cells $\pi^{\star} = ((a,b,c),(d),(e)) \in \mathcal{P}_5^\star.$
It corresponds to the Ferrers diagram shown in Figure~\ref{pi_in_Pstar}.

\begin{figure}[H]
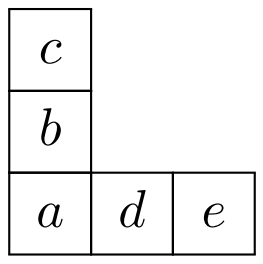

\centering
\[
\begin{ytableau}
c \\
b \\
a & d & e
\end{ytableau}
\]
\caption{Ferrers diagram of $\pi^{\star} = ((a,b,c),(d),(e))$.}
\label{pi_in_Pstar}
\end{figure}

Moreover, we have $\mathcal{L}_{\pi^\star}=\{a,b,c,d,e\}$.

\end{example}

For two partitions with distinguishable cells $\mu^\star$ and $\pi^\star$ with disjoint label sets, define their union $\mu^\star\cup\pi^\star$ by combining their columns and arranging
them in weakly decreasing order of length. The order of the cells
within each column is preserved. Among columns of the same length,
the original left-to-right order within each partition is preserved,
and all columns from $\pi^\star$ are placed before those from
$\mu^\star$. For example,
\[
\bigl((a,b,c),(d),(e)\bigr) \cup \bigl((A,B,C),(D,E),(F)\bigr) = \bigl((A,B,C),(a,b,c),(D,E),(F),(d),(e)\bigr).
\]

We now introduce four actions, $T_i$ for $i = 1,2,3,4$, acting on $\mathcal{P}_n^\star$. Each action is defined by choosing a rule for removing one cell from each column, followed by a rule for arranging the removed cells into a new column. There are two natural ways to remove cells: either extract the bottom cell $\pi_{i1}$ from each column or the top cell $\pi_{i n_i}$. Likewise, there are two natural methods for placing these cells: stack them in the new column from bottom to top, following their original left-to-right order, or arrange them from top to bottom. The four actions arise from all possible combinations of these removal and placement choices, and their formal definitions are given below.

\begin{definition}\label{def:actions}
Let $\pi \in \mathcal{P}_n^\star$ be
\[ \pi = \bigl((\pi_{11}, \pi_{12},  \dots, \pi_{1n_1}), (\pi_{21}, \pi_{22}, \dots, \pi_{2n_2}),\dots, (\pi_{k1}, \pi_{k2}, \dots, \pi_{kn_k})\bigr). \]
The actions $T_i:\mathcal{P}_n^\star \to \mathcal{P}_n^\star$, for $i=1,2,3,4$, are defined by
\begin{align*}
T_1(\pi) &= \bigl((\pi_{12},\dots,\pi_{1n_1}),(\pi_{22},\dots,\pi_{2n_2}), \dots,(\pi_{k2},\dots,\pi_{kn_k})\bigr) \;\cup\; (\pi_{11},\pi_{21},\dots,\pi_{k1}), \\
T_2(\pi) &= \bigl((\pi_{12},\dots,\pi_{1n_1}), (\pi_{22},\dots,\pi_{2n_2}), \dots,(\pi_{k2},\dots,\pi_{kn_k})\bigr) \;\cup\; (\pi_{k1},\pi_{k-1,1},\dots,\pi_{11}), 
\\
T_3(\pi) &= \bigl((\pi_{11},\dots,\pi_{1,n_1-1}), (\pi_{21},\dots,\pi_{2,n_{2}-1)}, \dots,(\pi_{k1},\dots,\pi_{k,n_k-1})\bigr) \;\cup\; (\pi_{1 n_1},\pi_{2 n_2},\dots,\pi_{k n_k}), \\
T_4(\pi) &= \bigl((\pi_{11},\dots,\pi_{1,n_1-1}), (\pi_{21},\dots,\pi_{2,n_2-1}), \dots,(\pi_{k1},\dots,\pi_{k,n_k-1})\bigr) \;\cup\; (\pi_{k n_k},\pi_{k-1,n_{k-1}},\dots,\pi_{1 n_1}).
\end{align*}

Empty columns are omitted. In these formulas, the tuple to the right
of $\cup$ is regarded as a one-column partition.
\end{definition}

In other words, $T_1$ and $T_2$ remove the bottom cell $\pi_{i1}$ from each column. They differ in how the removed cells are stacked: $T_1$ preserves the left-to-right order, while $T_2$ reverses it. Similarly, $T_3$ and $T_4$ remove the top cell $\pi_{i n_i}$ from each column; $T_3$ stacks them preserving the left-to-right order, while $T_4$ reverses it.

As proved in Proposition~3.4, a partition with distinguishable cells
is periodic under $T_i$ if and only if its underlying partition is
periodic under $T$. As in the classical Bulgarian solitaire case, for periodic partitions the new column is always at least as long as the remaining columns, and hence the union places it as the first column after applying each $T_i$. Figures~\ref{O_{T_1}} and~\ref{O_{T_2}} illustrate
the orbits of $\pi=((a,b,c),(d),(e))$ under $T_1$ and $T_2$, respectively. The other two actions are related to these by the reversal symmetry established in Proposition~\ref{prop:rev-relation}.

\begin{figure}[H]
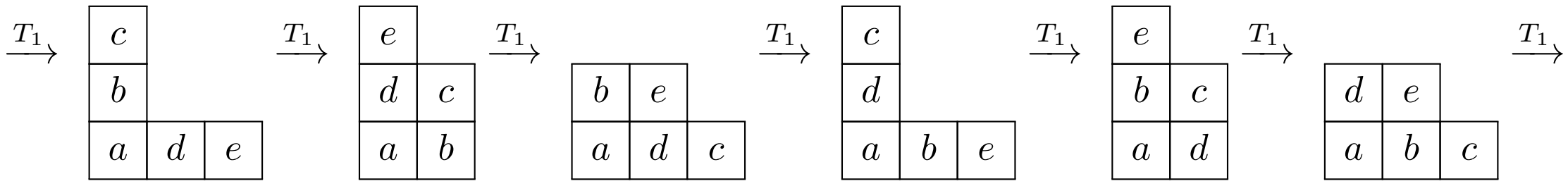


\[ \xrightarrow{T_1}\,\,\, \begin{ytableau} 
c  \\ 
b \\ 
a & d & e 
\end{ytableau} 
\xrightarrow{T_1}\,\,\, \begin{ytableau} 
e  \\ 
d & c\\ 
a & b 
\end{ytableau}
\xrightarrow{T_1}\,\,\, \begin{ytableau} 
\none & \none & \none \\
b & e \\
a & d & c 
\end{ytableau}
\xrightarrow{T_1}\,\,\, \begin{ytableau}
c  \\
d \\ 
a & b & e \end{ytableau}
\xrightarrow{T_1}\,\,\,\begin{ytableau} 
e  \\ 
b & c\\ 
a & d 
\end{ytableau}
\xrightarrow{T_1}\,\,\, \begin{ytableau} 
\none & \none & \none \\
d & e \\
a & b & c 
\end{ytableau}
\xrightarrow{T_1}\,\,\, \]

\caption{One orbit under the action $T_1$.}
\label{O_{T_1}}

\end{figure}

\begin{figure}[H]
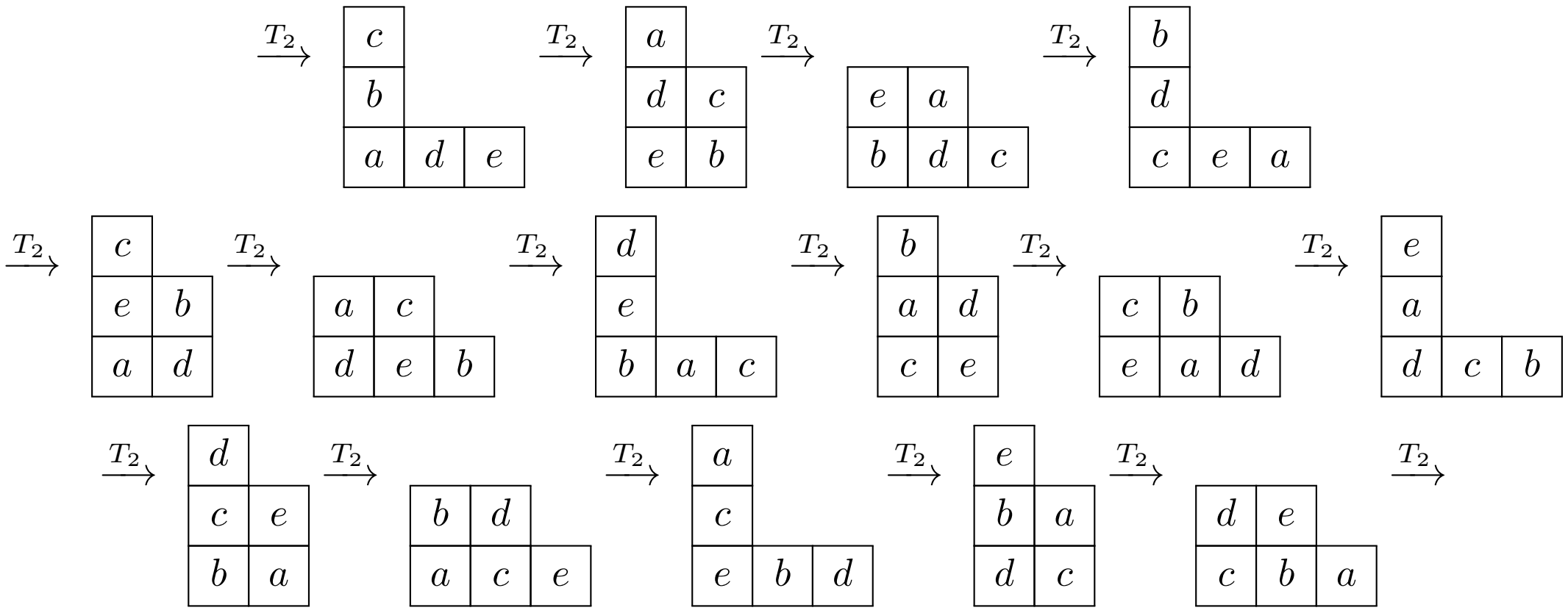


\[\xrightarrow{T_2}\,\,\,\begin{ytableau} 
c  \\ 
b \\ 
a & d & e 
\end{ytableau} 
\xrightarrow{T_2}\,\,\, \begin{ytableau} 
a  \\ 
d & c\\ 
e & b 
\end{ytableau}
\xrightarrow{T_2}\,\,\, \begin{ytableau} 
\none & \none & \none \\
e & a \\
b & d & c 
\end{ytableau}
\xrightarrow{T_2}\,\,\, \begin{ytableau}
b  \\
d \\ 
c & e & a \end{ytableau} \]
\[
\xrightarrow{T_2}\,\,\,\begin{ytableau} 
c  \\ 
e & b\\ 
a & d 
\end{ytableau}
\xrightarrow{T_2}\,\,\, \begin{ytableau} 
\none & \none & \none \\
a & c \\
d & e & b 
\end{ytableau}
\xrightarrow{T_2}\,\,\, \begin{ytableau} 
d  \\
e \\ 
b & a & c \end{ytableau}
\xrightarrow{T_2}\,\,\, \begin{ytableau} 
b  \\ 
a & d\\ 
c & e 
\end{ytableau}
\xrightarrow{T_2}\,\,\, \begin{ytableau} 
\none & \none & \none \\
c & b \\
e & a & d 
\end{ytableau}
\xrightarrow{T_2}\,\,\, \begin{ytableau}
e  \\
a \\ 
d & c & b \end{ytableau}\]
\[
\xrightarrow{T_2}\,\,\, \begin{ytableau} 
d  \\ 
c & e\\ 
b & a 
\end{ytableau}
\xrightarrow{T_2}\,\,\, \begin{ytableau} 
\none & \none & \none \\
b & d \\
a & c & e 
\end{ytableau}
\xrightarrow{T_2}\,\,\, \begin{ytableau}
a  \\
c \\ 
e & b & d \end{ytableau}
\xrightarrow{T_2}\,\,\, \begin{ytableau} 
e  \\ 
b & a\\ 
d & c 
\end{ytableau}
\xrightarrow{T_2}\,\,\, \begin{ytableau} 
\none & \none & \none \\
d & e \\
c & b & a 
\end{ytableau}
\xrightarrow{T_2}\,\,\, \]

\caption{One orbit under the action $T_2$.}
\label{O_{T_2}}
\end{figure}

\subsection{Connection to classical Bulgarian solitaire}\label{sec:connection-to-classical}

Each of the four actions $T_1,T_2,T_3,T_4$ removes one cell from every column and forms a new column of size equal to the number of columns.  
Consequently, although the labels are permuted differently, the multiset of column sizes evolves exactly as in the classical game.  
To make this precise, define the \emph{shape} map
\[
S:\mathcal{P}_n^\star \longrightarrow \mathcal{P}_n ,\qquad
S(\pi^\star)=\text{nonincreasing list of column lengths of }\pi^\star \in \mathcal{P}_n^\star  .
\]

\begin{proposition}\label{prop:shape-intertwines-all}
For every $\pi^\star\in\mathcal{P}_n^\star$ and every $i\in\{1,2,3,4\}$,
\[
S\bigl(T_i(\pi^\star)\bigr)=T\bigl(S(\pi^\star)\bigr),
\]
where $T:\mathcal{P}_n\to\mathcal{P}_n$ is the classical Bulgarian solitaire map.
In other words, the following diagram commutes:
\[
\begin{tikzcd}
\mathcal{P}_n^\star \ar[r,"T_i"] \ar[d,"S"'] & \mathcal{P}_n^\star \ar[d,"S"] \\
\mathcal{P}_n \ar[r,"T"] & \mathcal{P}_n
\end{tikzcd}
\]
\end{proposition}

\begin{proof}
Let $S(\pi^\star)=(a_1,\ldots,a_k)$. Each $T_i$ removes one cell
from every column and forms a new column from the $k$ removed
cells. Thus
\[
S(T_i(\pi^\star))
=(a_1-1,\ldots,a_k-1)\cup(k)
=T(S(\pi^\star)),
\]
where zero parts are omitted and the remaining parts are arranged
in weakly decreasing order.
\end{proof}

\begin{proposition}\label{prop:periodic shape} For each $i\in\{1,2,3,4\}$, a partition with distinguishable cells
$\pi^\star$ is periodic under $T_i$ if and only if its shape
$S(\pi^\star)$ is periodic under $T$.
\end{proposition}

\begin{proof}
\noindent
\((\Rightarrow)\) First suppose that \(\pi^\star\) is periodic under \(T_i\). Then, for some \(\ell>0\),
\[
T_i^\ell(\pi^\star)=\pi^\star .
\]
Using Proposition~\ref{prop:shape-intertwines-all}, namely \(S\circ T_i=T\circ S\), we get
\[
T^\ell(S(\pi^\star))
=
S(T_i^\ell(\pi^\star))
=
S(\pi^\star).
\]
Thus \(S(\pi^\star)\) is periodic under \(T\).

\medskip
\noindent
\((\Leftarrow)\)
Conversely, suppose that \(\lambda = S(\pi^\star)\) is periodic under \(T\). Then there exists \(p>0\) such that \(T^{p}(\lambda) = \lambda\). For any partition \(\nu\in\mathcal{P}_n\), denote the fiber over \(\nu\) by
\[
F_\nu = \{\mu^\star\in\mathcal{P}_n^\star \mid S(\mu^\star) = \nu\}.
\]
Each \(F_\nu\) is finite; indeed \(|F_\nu| = n!\) because it consists of all bijections from a set of \(n\) distinguishable labels to the cells of the Ferrers diagram of \(\nu\).

By Proposition~\ref{prop:shape-intertwines-all}, we have \(T_i(F_\nu) \subseteq F_{T(\nu)}\). Moreover, the map
\[
T_i : F_\nu \longrightarrow F_{T(\nu)}
\]
is injective: the source shape $\nu$ is fixed, so the newly created column has known length $\ell(\nu)$. By the tie-breaking rule, it is the leftmost column of that length. Removing it leaves the nonempty remnants of the original columns in their original order. Reading its labels in the order prescribed by $T_i$ and reinserting them reconstructs the original configuration, including any singleton columns that disappeared. Since \(|F_\nu| = |F_{T(\nu)}| = n!\), the injection is a bijection.

Consequently, \(T_i^{p}\) maps the fiber \(F_\lambda\) bijectively onto \(F_{T^{p}(\lambda)} = F_\lambda\). Thus \(T_i^{p}\) is a permutation of the finite set \(F_\lambda\). Any permutation of a finite set has finite order, so there exists \(m\ge 1\) such that \((T_i^{p})^{m}(\pi^\star) = \pi^\star\). Therefore \(\pi^\star\) is periodic under \(T_i\).

\end{proof}

The following is an immediate consequence of
Proposition~\ref{prop:periodic shape}.

\begin{corollary}\label{prop:preperiod}
For $\pi^\star\in\mathcal{P}_n^\star$ and
$i\in\{1,2,3,4\}$, let $k_i$ and $k$ be the least nonnegative
integers such that $T_i^{k_i}(\pi^\star)$ and
$T^k(S(\pi^\star))$ are periodic, respectively. Then $k_i=k$.
In particular, the maximum preperiod under each of the four
actions is $D(n)$.
\end{corollary}

\subsection{\texorpdfstring{Connections between $T_1, T_2, T_3$ and $T_4$} {Connections between T1, T2, T3 and T4}} \label{sec:orbits}

The four maps introduced above are related pairwise through a natural reversal operation, which we now define.  

\begin{definition}\label{def:rev}
For $\pi=(\pi_1,\dots,\pi_k)\in\mathcal{P}_n^\star$, where $\pi_i=(\pi_{i1},\dots,\pi_{in_i})$, define
\[
\rev:\mathcal{P}_n^\star\to\mathcal{P}_n^\star,\qquad
\rev(\pi)=\bigl(\rev(\pi_1),\dots,\rev(\pi_k)\bigr),
\]
where $\rev(\pi_i)=(\pi_{in_i},\pi_{in_i-1},\dots,\pi_{i1})$ for $1\le i\le k$.
\end{definition}

\begin{proposition}\label{prop:rev-relation}
Let $\pi \in \mathcal{P}^\star_n$. The following diagrams commute:
\[
\begin{tikzcd}
\pi \ar[r,"T_1"] \ar[d,"\rev"'] & T_1(\pi) \ar[d,"\rev"] \\
\rev(\pi) \ar[r,"T_4"] & T_4(\rev(\pi))
\end{tikzcd}
\qquad\qquad
\begin{tikzcd}
\pi \ar[r,"T_2"] \ar[d,"\rev"'] & T_2(\pi) \ar[d,"\rev"] \\
\rev(\pi) \ar[r,"T_3"] & T_3(\rev(\pi))
\end{tikzcd}
\]
    i.e.,  $\rev \circ T_1 = T_4 \circ \rev$  and  $\rev \circ T_2 = T_3 \circ \rev$.
   
\end{proposition}

\begin{proof}
Reversing a column after removing its bottom cell gives the same
column as removing its top cell after reversal. Hence the columns
remaining from $\pi$ are identical in
$\operatorname{rev}(T_1(\pi))$ and
$T_4(\operatorname{rev}(\pi))$. The new columns also agree:
both contain the original bottom cells in right-to-left order
when read from bottom to top. Since reversal preserves column
lengths and the left-to-right order of columns, the sorting and
tie-breaking rules give the same arrangement in both constructions.
Thus $\operatorname{rev}\circ T_1=T_4\circ\operatorname{rev}$.

The same argument applies to $T_2$ and $T_3$, except that the
original bottom cells appear in the new column in left-to-right
order when read from bottom to top. Therefore
$\operatorname{rev}\circ T_2=T_3\circ\operatorname{rev}$.
\end{proof}

\begin{corollary}\label{same_orbit_sizes} For every periodic $\pi\in\mathcal{P}^\star_n$,
    $|O_{T_1}(\pi)|= |O_{T_4}(\pi)|$  and   $|O_{T_2}(\pi)|= |O_{T_3}(\pi)|$.
\end{corollary}

\begin{proof}
By Proposition~3.7, $\operatorname{rev}$ sends the $T_1$-orbit of
$\pi$ bijectively to the $T_4$-orbit of $\operatorname{rev}(\pi)$,
and similarly sends the $T_2$-orbit of $\pi$ bijectively to the
$T_3$-orbit of $\operatorname{rev}(\pi)$.

For each $i$, the map $T_i$ commutes with every permutation of the
labels, since its rules depend only on the positions of the cells.
Consequently, any two labelings of the same shape have the same
orbit size under $T_i$. Since $\pi$ and $\operatorname{rev}(\pi)$
have the same shape, the stated equalities follow.
\end{proof}
Therefore, it suffices to study the orbit sizes of periodic $\pi\in\mathcal{P}^\star_n$ under the actions $T_1$ and $T_2$.

\section{Initial analysis of the orbit sizes} \label{sec:orbit sizes}

In this section, we first determine the orbit size under $T_1$, then focus on $T_2$ by developing a coloring algorithm on the Ferrers diagram. 
The coloring algorithm reveals the permutation $\sigma_\pi$ induced by $T_2^{m+1}$ on the labeled cells, and we prove that the colored regions correspond exactly to the cycles of $\sigma_\pi$. 
We then compute the order of $\sigma_\pi$ explicitly in terms of the parameters $m$ and $m'$, obtaining a bound on the orbit size under $T_2$, with equality when the underlying shape has period $m+1$. The section concludes with several illustrative examples.

\subsection{\texorpdfstring
  {The orbit size of $\pi$ under the action $T_1$}
  {The orbit size of pi under the action T1}} \label{subsec:T1}

\begin{theorem}\label{thm:orbit size T1}
Let $n$ and $m$ be positive integers such that $\Delta_m < n \leq \Delta_{m+1}$. Then every periodic orbit of $\mathcal{P}_n^\star$ under $T_1$ has size $\lcm(1,2,\ldots,m+1)$.
\end{theorem}

\begin{proof} Suppose $\pi \in \mathcal{P}_n^\star$ is a periodic partition with distinguishable cells. By Theorem~\ref{periodic points under T}, $\pi$ contains cells in $\diag[k]$ for $1 \leq k \leq m+1$. We track how these cells move under $T_1$. The cell in position $(i,j)$ moves to position $(i+1,j-1)$, for $j>1$. The cell in position $(i,1)$ moves to position $(1,i)$. 

These moves cyclically permute the $k$ positions of $\operatorname{diag}[k]$,
so each cell on this diagonal first returns to its original position
after exactly $k$ iterations. Since every diagonal
$\operatorname{diag}[k]$ with $1\leq k\leq m+1$ contains at least one
cell, and all labels are distinct, the orbit size is
$\operatorname{lcm}(1,2,\ldots,m+1)$.
\end{proof}

\subsection{\texorpdfstring
  {The orbit size of $\pi$ under the action $T_2$ }{The orbit size of pi under the action T2 }} \label{subsec:T2}

To better describe the movement of individual cells under the action of $T_2$, we introduce coordinates. For a partition $\pi$, let $\pi(i,j)$ denote the cell in the $i^{\text{th}}$ column and $j^{\text{th}}$ row of its Ferrers diagram. The following proposition describes the new location of this cell after applying $T_2$.

\begin{proposition}\label{actionT_2inverse}
Let $n$ and $m$ be positive integers such that $\Delta_m \leq n < \Delta_{m+1}$ and let $\pi \in \mathcal{P}_n^\star$ be a periodic point. Then
\[
\pi(i,j)=
\begin{cases} 
(T_2\pi)(i+1,j-1) & \text{if } j>1, \\[4pt]
(T_2\pi)(1,m+1-i) & \text{if } j=1 \text{ and } \pi \text{ has } m \text{ cells in its bottom row}, \\[4pt]
(T_2\pi)(1,m+2-i) & \text{if } j=1 \text{ and } \pi \text{ has } m+1 \text{ cells in its bottom row}.
\end{cases}
\]
\end{proposition}

\subsubsection{\texorpdfstring
  {The permutation $\sigma_\pi$ }{The induced permutation}}

Let $\pi\in\mathcal{P}_n^\star$ be periodic. Since the period of
$S(\pi)$ under $T$ divides $m+1$, the configurations $\pi$ and
$T_2^{m+1}(\pi)$ have the same shape. Comparing the labels in
corresponding positions therefore defines a permutation of
$\mathcal{L}_\pi$.

\begin{definition}\label{actionT_2inverse=sigmapi}
Let $\pi \in \mathcal{P}_n^\star$ be a periodic partition with distinguishable cells, and let $m \in \mathbb{N}$ be such that $\Delta_m \leq n <  \Delta_{m+1}$. Define $\sigma_\pi: \mathcal{L}_\pi \to \mathcal{L}_\pi$ by
\[
\sigma_\pi(\pi(i,j)) = (T_2^{m+1}\pi)(i,j).
\]
\end{definition}

Thus, $\sigma_\pi(a)=b$ means that after $m+1$ iterations of $T_2$, the
label $b$ occupies the position initially occupied by $a$.

\begin{example}\label{ex:five-cell-permutation}
For $\pi=((a,b,c),(d),(e))$ in Figure~\ref{O_{T_2}}, we have $m=2$ and
$T_2^3(\pi)=((c,d,b),(e),(a))$. Comparing the labels in
corresponding positions gives $\sigma_\pi=(a\,c\,b\,d\,e)$.
\end{example}

Recall that the order $\operatorname{ord}(\sigma)$ of a permutation $\sigma \in S_n$ is the least positive integer $k$ such that $\sigma^k = \mathrm{id}$; equivalently, it is the least common multiple of the lengths of the cycles in its disjoint cycle decomposition.

The orbit of a partition $\pi \in \mathcal{P}_n^\star$ under $T_2$ exhibits a dual periodicity: the Ferrers diagram returns to its original shape every $m+1$ iterations of $T_2$, while the permutation $\sigma_\pi$ governs the cyclic rearrangement of the labels within this fixed shape. We therefore obtain the following result.

\begin{theorem}\label{thm:maximal orbit size}
Let $\pi \in \mathcal{P}_n^\star$ be a periodic partition under $T_2$, and let $m \in \mathbb{N}$ be such that  $\Delta_m \leq n < \Delta_{m+1}$. Then
\[
|O_{T_2}(\pi)| \bigm| (m+1)\operatorname{ord}(\sigma_\pi).
\] Moreover, if the period of the underlying shape $S(\pi)$ under $T$ is $\ell = m+1$, then 
\[
|O_{T_2}(\pi)| =   (m+1)\operatorname{ord}(\sigma_\pi).
\]
\end{theorem}

\begin{proof}
Because the rules of $T_2$ depend only on positions, successive
applications of $T_2^{m+1}$ act on the labels by successive powers
of $\sigma_\pi$. Consequently,
$T_2^{(m+1)\operatorname{ord}(\sigma_\pi)}(\pi)=\pi$, so the period
of $\pi$ divides $(m+1)\operatorname{ord}(\sigma_\pi)$.

If the underlying shape has period $m+1$, every return time of
$\pi$ is a multiple of $m+1$. Moreover,
$T_2^{k(m+1)}(\pi)=\pi$ if and only if $\sigma_\pi^k$ is the
identity. The least such positive $k$ is
$\operatorname{ord}(\sigma_\pi)$, proving the equality.
\end{proof}

For the configuration in Example~\ref{ex:five-cell-permutation}, the underlying shape has period $m+1=3$ and $\operatorname{ord}(\sigma_\pi)=5$. Theorem~\ref{thm:maximal orbit size} therefore gives $|O_{T_2}(\pi)|=3\cdot5=15$, in agreement with
Figure~\ref{O_{T_2}}.

\subsubsection{Coloring algorithm}\label{subsec:algo}

The goal of the following subsection will be to describe the order of $\sigma_\pi$. We begin by introducing notation that will be used throughout.
\begin{definition}
Let $\pi \in \mathcal{P}_n^\star$ be a periodic partition, and let $m \in \mathbb{N}$ be such that $\Delta_m \leq n < \Delta_{m+1}$. The cells lying on the $(m+1)$-th diagonal, denoted $\diag[m+1]$, are called \emph{extra cells}. Their number is
\[
m' = n - \frac{m(m+1)}{2},
\]
where $0 \leq m' \leq m$.

 A row (respectively, column) containing an extra cell is called a \emph{long row} (respectively, \emph{long column}), while a row (respectively, column) containing no extra cell is called a \emph{short row} (respectively, \emph{short column}).
\end{definition}

We show that the cycle structure of \(\sigma_\pi\) is encoded by a coloring of the cells of \(\pi\). In the following description of the coloring procedure, we assume \(m' > 0\). When \(m' = 0\), the diagonal \(\diag[m+1]\) contains no extra cells, and every cell remains uncolored.

Let
\[
E = \big( (i_1,h_1), \ldots, (i_{m'},h_{m'}) \big)
\]
be the ordered list of extra cells of $\pi$ along the diagonal $\diag[m+1]$, with $i_1 < i_2 < \cdots < i_{m'}$.

Define a map $X : \mathbb{N}^2 \times \mathbb{N}^2 \to \mathbb{N}^2$ by
\[
X\big((a_1,b_1),(a_2,b_2)\big) = (a_1,b_2),
\]
which returns the intersection of the column of $(a_1,b_1)$ and the row of $(a_2,b_2)$.

For $0 \leq k \leq m'-1$, define $X^k(E)$ inductively by $X^0(E) = E$ and
\[
X^k(E) = X\big(X^{k-1}(E)\big),
\]
where $X$ is applied to consecutive pairs of entries in the input ordered list. That is, if
\[
X^{k-1}(E) = \big( (c_1,d_1), \ldots, (c_r,d_r) \big),
\]
then
\[
X^k(E) = \big( X((c_1,d_1),(c_2,d_2)), \ldots, X((c_{r-1},d_{r-1}),(c_r,d_r)) \big).
\] Note that $X^{m'-1}(E) = \big( (i_1,h_{m'}) \big)$.

We now define a coloring of the cells of $\pi$. Let $\{C_1, \dots, C_{\lceil m'/2 \rceil}\}$ be a set of distinct colors. For each $1 \leq r \leq \lceil m'/2 \rceil$, we perform an \emph{initial coloring} by assigning color $C_r$ to the cells in
\[
X^{r-1}(E) \cup X^{m'-r}(E).
\]

After all initial colors have been assigned, extend the color of
each initially colored cell straight downward in its column and
straight leftward in its row, stopping before the first previously
colored cell or at the boundary of the diagram. Newly colored cells
do not initiate further \emph{propagation}.

We also use $C_i$ to denote the color class consisting of all cells
assigned color $C_i$ by this procedure.

\begin{example}\label{ex:pi31}
Consider the colored partition in Figure~\ref{fig:colored pi31} with $n=31$, $m=7$, and $m'=3$. The diagonal $\diag[8]$ contains the extra cells $(1,8)$, $(4,5)$, and $(7,2)$.

We have
\[
E = \big((1,8),(4,5),(7,2)\big), \quad 
X^1(E) = \big((1,5),(4,2)\big), \quad 
X^2(E) = \big((1,2)\big).
\]
Since $m'=3$, we use $\lceil m'/2 \rceil = 2$ colors. The initial coloring is given by
\[
C_1 (yellow): X^0(E) \cup X^2(E) = \{(1,8),(4,5),(7,2),(1,2)\}, \quad
C_2 (blue): X^1(E) = \{(1,5),(4,2)\}.
\]
The coloring is then extended downward in columns and leftward in rows, as shown in Figure~\ref{fig:colored pi31}.

\begin{figure}[H]
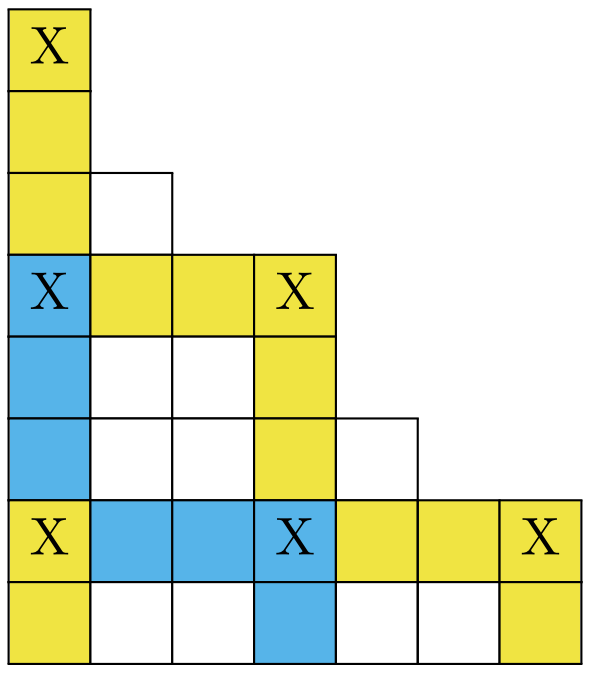
   

\[
\begin{ytableau}
*(yellow) \text{X} \\  
*(yellow) \\     
*(yellow) & *(white) \\ 
*(green)\text{X} & *(yellow) &*(yellow) &*(yellow)\text{X} \\
*(green) & *(white) &*(white) &*(yellow) \\
*(green) & *(white) &*(white) &*(yellow)&*(white) \\
 *(yellow)\text{X} & *(green) &*(green) &*(green)\text{X} & *(yellow) &*(yellow) &*(yellow)\text{X}\\
*(yellow)& *(white) &*(white) & *(green) &*(white) &*(white) & *(yellow)  \\                      
\end{ytableau}
\]
   \caption{The cell coloring procedure applied to a $\pi \in \mathcal{P}^\star_{31}$. The cells marked with `X' indicate the initially colored intersection points before color propagation.}
   \label{fig:colored pi31}
\end{figure}
\end{example}

\begin{figure}[H]
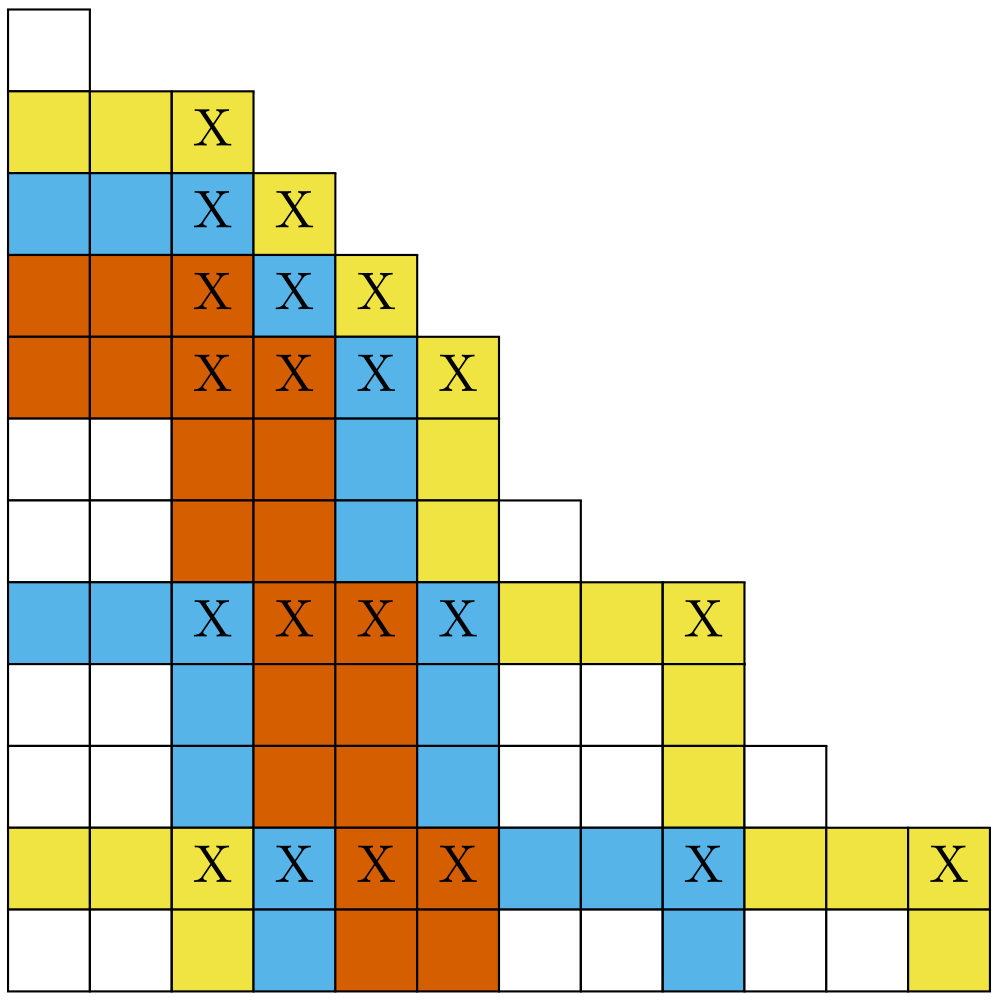

\[
\begin{ytableau}
*(white) \\  
*(yellow) & *(yellow) & *(yellow) \text{X} \\     
*(green) & *(green) & *(green) \text{X} & *(yellow) \text{X} \\ 
*(orange) & *(orange) & *(orange) \text{X} & *(green) \text{X} & *(yellow) \text{X} \\
*(orange) & *(orange) & *(orange) \text{X} & *(orange) \text{X} & *(green) \text{X} & *(yellow) \text{X} \\
*(white) & *(white) & *(orange) & *(orange) & *(green) & *(yellow) \\
*(white) & *(white) & *(orange) & *(orange) & *(green) & *(yellow) & *(white) \\
*(green) & *(green) & *(green) \text{X} & *(orange) \text{X} & *(orange) \text{X} & *(green) \text{X} & *(yellow) & *(yellow) & *(yellow) \text{X} \\
*(white) & *(white) & *(green) & *(orange) & *(orange) & *(green) & *(white) & *(white) & *(yellow) \\
*(white) & *(white) & *(green) & *(orange) & *(orange) & *(green) & *(white) & *(white) & *(yellow) & *(white) \\
*(yellow) & *(yellow) & *(yellow) \text{X} & *(green) \text{X} & *(orange) \text{X} & *(orange) \text{X} & *(green) & *(green) & *(green) \text{X} & *(yellow) & *(yellow) & *(yellow) \text{X} \\
*(white) & *(white) & *(yellow) & *(green) & *(orange) & *(orange) & *(white) & *(white) & *(green) & *(white) & *(white) & *(yellow) \\
\end{ytableau}
\]

\caption{The cell coloring procedure applied to a $\pi \in \mathcal{P}^\star_{84}$. The cells marked with `X' indicate the initially colored intersection points before color propagation.}
\label{fig:colored pi86}
\end{figure}

\begin{example}\label{ex:pi86}
For $n=84$, $m=12$, and $m'=6$, consider the periodic partition
whose extra cells are
\[
E=((3,11),(4,10),(5,9),(6,8),(9,5),(12,2)).
\]
The initial coloring assigns $C_1$ (yellow) to
$E\cup X^5(E)$, $C_2$ (blue) to $X^1(E)\cup X^4(E)$,
and $C_3$ (orange) to $X^2(E)\cup X^3(E)$.
The resulting coloring after propagation is shown in
Figure~\ref{fig:colored pi86}.
\end{example}

We finish this subsection with a couple of observations about the resulting coloring which will be useful in the forthcoming analyses. 

First, some cells of $\pi$ may remain uncolored by the coloring algorithm. In fact, a cell remains uncolored if and only if it lies in a short row and a short column of $\pi$. Indeed, every long row and every long column contains an extra cell, from which the coloring propagates leftward or downward,
while a cell lying in both a short row and a short column cannot be reached by either type of propagation.The total number of uncolored cells depends only on $m$ and $m'$ as given in the following lemma.

\begin{lemma}\label{lem:uncolored-count}
Let \(\pi \in \mathcal{P}_n^\star\) be a periodic point with \(n = \Delta_m + m'\) where \(0 \leq m' \leq m\). Then the number of uncolored cells in \(\pi\) under the coloring algorithm is exactly the triangular number
\[
\Delta_{m-m'} = \frac{(m-m')(m-m'+1)}{2}.
\]
\end{lemma}

\begin{proof}
There are $m+1-m'$ unoccupied positions on
$\operatorname{diag}[m+1]$. Choose two such positions with column
indices $a<b$. The cell $(a,m+2-b)$ lies in a short column and a
short row, and it belongs to the Ferrers diagram because
$a+(m+2-b)\leq m+1$.

Conversely, every uncolored cell $(p,q)$ lies in a short column
and a short row, so it determines two unoccupied positions on
$\operatorname{diag}[m+1]$ with column indices
$a=p$ and $b=m+2-q$. Since $(p,q)$ is not an extra cell,
$p+q\leq m+1$, and hence $a<b$.

Thus uncolored cells are in bijection with pairs of unoccupied
positions on $\operatorname{diag}[m+1]$, and their number is
$\binom{m+1-m'}{2}=\Delta_{m-m'}$.
\end{proof}

Secondly, the initial coloring pairs the sets $X^t(E)$ and $X^{m'-1-t}(E)$ for $t=0,\ldots,\lceil m'/2\rceil-1$, from which
the colors propagate. Each of these sets, together with the cells colored by propagation from it, forms a chain in the weak southeast
order. When $m'$ is even, each color class consists of two distinct chains (see Figure~\ref{fig:colored pi86}). When $m'$ is odd, each color class likewise consists of two distinct chains, except for the last color class, corresponding to $t=\lfloor m'/2\rfloor$, for which the two chains coincide because $X^t(E)=X^{m'-1-t}(E)$. This occurs in Figure~\ref{fig:colored pi31}, where $m'=3$ and the last color class consists of a single chain.

\subsubsection{\texorpdfstring
  {The permutation $\sigma_\pi^{\mathrm{col}}$}{The colored induced}}

We now relate the color classes of $\pi$ to the cycle structure of $\sigma_\pi$ induced by the action of $T_2$. The coloring not only partitions the Ferrers diagram but also encodes the cycle decomposition of $\sigma_\pi$: each color class corresponds to one cycle, while the uncolored labels are fixed. To make this precise, we define a permutation $\sigma_\pi^{\mathrm{col}}$ determined by the coloring and show that $\sigma_\pi^{\mathrm{col}} = \sigma_\pi$.

We define the permutation $\sigma_\pi^{\mathrm{col}}$ by specifying its disjoint cycle decomposition. For each color $C_i$, $1 \leq i \leq \lceil m'/2 \rceil$, we construct a cycle $\sigma_\pi^{C_i}$ as follows.

Using the heights $h_1,\ldots,h_{m'}$ from $E$, set $l_i=h_{m'-i+1}$. To construct $\sigma_\pi^{C_i}$, first list the
labels in the chain generated by $X^{i-1}(E)$ in weak southeast order, starting at $(1,h_i)$ and ending at $(m+2-l_i,1)$. When the two chains are distinct, append the labels in the chain generated by $X^{m'-i}(E)$ in the same order, starting at $(1,l_i)$ and ending at $(m+2-h_i,1)$. Read the resulting list
cyclically, so that the final label maps to the first.

Thus, each $\sigma_\pi^{C_i}$ is a cycle beginning at $(1,h_i)$ and ending at $(m+2-h_i,1)$.

\begin{definition}\label{sigmapicol}
Let $\pi\in\mathcal{P}_n^\star$ be a periodic partition with label
set $\mathcal{L}_\pi$, and color its cells using the algorithm in
Section~4.2.2. The colored permutation
$\sigma_\pi^{\mathrm{col}}\colon\mathcal{L}_\pi\to\mathcal{L}_\pi$
is defined by
\[
\sigma_\pi^{\mathrm{col}}
=
\sigma_\pi^{C_1}\sigma_\pi^{C_2}\cdots
\sigma_\pi^{C_{\lceil m'/2\rceil}},
\]
where the cycles $\sigma_\pi^{C_i}$ are as constructed above.
Every label of an uncolored cell is fixed.
\end{definition}

\begin{remark}\label{r1: moves of sigmapicol}
For each cell $(p,q)$, the following rules identify the label
mapped to $\pi(p,q)$ under $\sigma_\pi^{\mathrm{col}}$.

\begin{enumerate}
    \item \textbf{Fixed points}: If row $q$ and column $p$ are both short, then
    \[
    \sigma_\pi^{\mathrm{col}}(\pi(p,q)) = \pi(p,q).
    \]

    \item \textbf{Downward $\shortarrow{6}$ move}: If row $q$ is short and column $p$ is long, then
    \[
    \sigma_\pi^{\mathrm{col}}(\pi(p,q+1)) = \pi(p,q).
    \]

    \item \textbf{Diagonal $\shortarrow{7}$ move}: If row $q$ is long and column $p-1$ is long ($p \geq 2$), then
    \[
    \sigma_\pi^{\mathrm{col}}(\pi(p-1,q+1)) = \pi(p,q).
    \]

    \item \textbf{Rightward $\shortarrow{0}$ move}: If row $q$ is long and column $p-1$ is short ($p \geq 2$), then
    \[
    \sigma_\pi^{\mathrm{col}}(\pi(p-1,q)) = \pi(p,q).
    \]

    \item \textbf{Jump \(\circlearrowright\) move}: If \(p=1\) and row \(q\) is long, then \(q=l_i\) or \(q=h_i\), and
    \[
    \sigma_\pi^{\mathrm{col}}(\pi(m+2-l_i,1)) = \pi(1,l_i), \qquad
    \sigma_\pi^{\mathrm{col}}(\pi(m+2-h_i,1)) = \pi(1,h_i).
    \]
\end{enumerate}
\end{remark}

\begin{figure}[H]
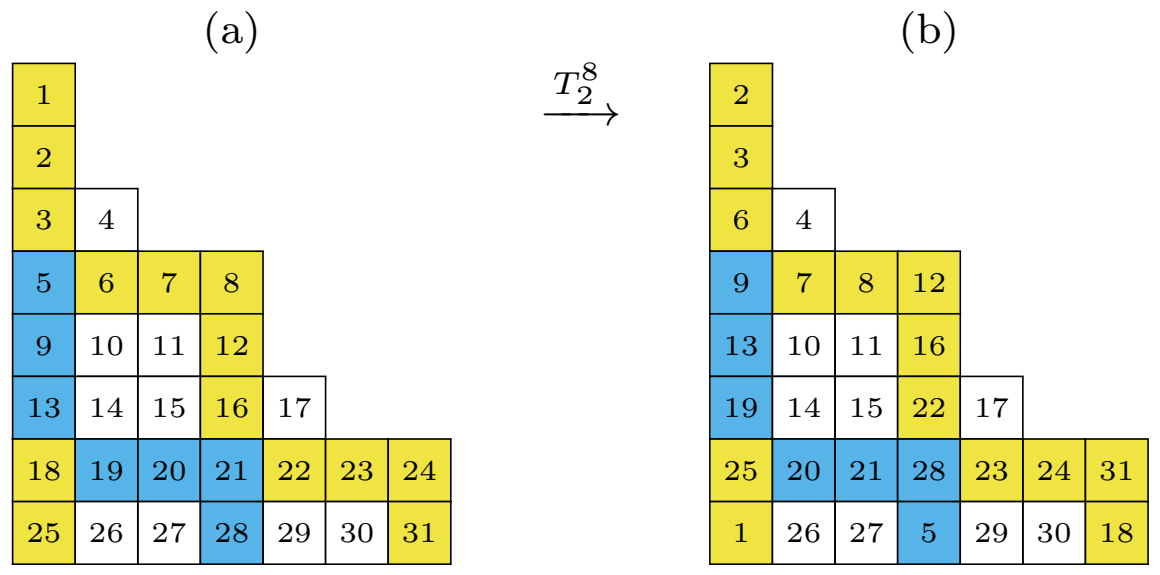

\centering
\footnotesize
\begin{tabular}{@{}c@{\hspace{5mm}}c@{\hspace{5mm}}c@{}}
(a) & & (b) \\ 

\scalebox{0.9}{\begin{ytableau}
*(yellow) {\scriptstyle 1} \\  
*(yellow) {\scriptstyle 2} \\     
*(yellow) {\scriptstyle 3} & *(white) {\scriptstyle 4} \\ 
*(green) {\scriptstyle 5} & *(yellow) {\scriptstyle 6} & *(yellow) {\scriptstyle 7} & *(yellow) {\scriptstyle 8} \\
*(green) {\scriptstyle 9} & *(white) {\scriptstyle 10} & *(white) {\scriptstyle 11} & *(yellow) {\scriptstyle 12} \\
*(green) {\scriptstyle 13} & *(white) {\scriptstyle 14} & *(white) {\scriptstyle 15} & *(yellow) {\scriptstyle 16} & *(white) {\scriptstyle 17} \\
 *(yellow) {\scriptstyle 18} & *(green) {\scriptstyle 19} & *(green) {\scriptstyle 20} & *(green) {\scriptstyle 21} & *(yellow) {\scriptstyle 22} & *(yellow) {\scriptstyle 23} & *(yellow) {\scriptstyle 24}\\
*(yellow) {\scriptstyle 25} & *(white) {\scriptstyle 26} & *(white) {\scriptstyle 27} & *(green) {\scriptstyle 28} & *(white) {\scriptstyle 29} & *(white) {\scriptstyle 30} & *(yellow) {\scriptstyle 31}                      
\end{ytableau}}
&$\;\xrightarrow{T_2^8}\;$ &

\scalebox{0.9}{\begin{ytableau}
*(yellow) {\scriptstyle 2} \\  
*(yellow) {\scriptstyle 3} \\     
*(yellow) {\scriptstyle 6} & *(white) {\scriptstyle 4} \\ 
*(green) {\scriptstyle 9} & *(yellow) {\scriptstyle 7} & *(yellow) {\scriptstyle 8} & *(yellow) {\scriptstyle 12} \\
*(green) {\scriptstyle 13} & *(white) {\scriptstyle 10} & *(white) {\scriptstyle 11} & *(yellow) {\scriptstyle 16} \\
*(green) {\scriptstyle 19} & *(white) {\scriptstyle 14} & *(white) {\scriptstyle 15} & *(yellow) {\scriptstyle 22} & *(white) {\scriptstyle 17} \\
 *(yellow) {\scriptstyle 25} & *(green) {\scriptstyle 20} & *(green) {\scriptstyle 21} & *(green) {\scriptstyle 28} & *(yellow) {\scriptstyle 23} & *(yellow) {\scriptstyle 24} & *(yellow) {\scriptstyle 31}\\
*(yellow) {\scriptstyle 1} & *(white) {\scriptstyle 26} & *(white) {\scriptstyle 27} & *(green) {\scriptstyle 5} & *(white) {\scriptstyle 29} & *(white) {\scriptstyle 30} & *(yellow) {\scriptstyle 18}                      
\end{ytableau}}
\end{tabular}
\caption{The configurations $\pi$ in (a) and $T_2^8(\pi)$ in (b), with cells colored by the coloring algorithm. Comparing labels
in corresponding positions gives the induced permutation $\sigma_\pi$.}
\label{fig:T2_iterations_9}
\end{figure}

\begin{example}
Consider the labeling in Figure~\ref{fig:T2_iterations_9}(a)
of the partition colored in Figure~\ref{fig:colored pi31},
where $n=31$, $m=7$, and $m'=3$.
The two yellow chains give a cycle of length $14$, while
the single blue chain gives a cycle of length $7$:
\[
\sigma_\pi^{\mathrm{col}}
=
(1\;2\;3\;6\;7\;8\;12\;16\;22\;23\;24\;31\;18\;25)\;(5\;9\;13\;19\;20\;21\;28).
\]
All uncolored labels are fixed.
Figure~\ref{fig:T2_iterations_9}(b) shows $T_2^8(\pi)$.
Comparing the labels in corresponding positions in the two
panels gives the same permutation, verifying that
$\sigma_\pi=\sigma_\pi^{\mathrm{col}}$ in this example.
\end{example}

\subsubsection{\texorpdfstring
  {$\sigma_{\pi}^{col} = \sigma_{\pi}$}{Equality of the coloring and induced permutations}}

We now prove that $\sigma_\pi^{\mathrm{col}} = \sigma_\pi$. This requires analyzing the behavior of cells under the action $T_2$ and verifying that the coloring and cycle-tracing procedures capture the induced dynamics. We begin with a lemma describing the movement of individual cells.

\begin{lemma}\label{cell movement}
Let $\pi \in \mathcal{P}_n^\star$ be a periodic partition, and let $a$ be a label in cell $(p,q)$ of $\pi$. Under repeated applications of $T_2$, the following hold:

\begin{enumerate}[label=(\arabic*)]
    \item If $(p,q)$ is in a short row and a short column, then the label $a$ remains in a short row and a short column.
    
    \item If $(p,q)$ is in a short row and a long column, then the label $a$ remains in a short row and a long column until it reaches the bottom row, after which it moves to a long row and a short column.
    
    \item If $(p,q)$ is in a long row and column $p-1$ is long ($p \geq 2$), then the label $a$ remains in a long row until it reaches the bottom row, after which it moves to a long row and a long column.
    
    \item If $(p,q)$ is in a long row and column $p-1$ is short ($p \geq 2$), then the label $a$ remains in a long row until it reaches the bottom row, after which it moves to a short row and a long column.
\end{enumerate}
\end{lemma}

\begin{proof}
We analyze each case using Proposition~\ref{actionT_2inverse} and the structure of periodic partitions.

Recall that a row (respectively, column) is short if it does not contain an extra cell, equivalently if there is an empty position on the diagonal $\diag[m+1]$ in that row (respectively, column).

\medskip

\noindent\textbf{(1)} Suppose $(p,q)$ lies in a short row and a short column. Then there are empty positions at $(m+2-q,q)$ and $(p,m+2-p)$. 

For $q>1$, the map $T_2$ moves the label $a$ from $(p,q)$ to $(p+1,q-1)$, while these empty positions move to $(m+3-q,q-1)$ and $(p+1,m+1-p)$, respectively. Thus $a$ remains in a row and column that are both short.

If $q=1$, the label $a$ moves to $(1,m+1-p)$. The corresponding empty positions ensure that this row and column remain short, so the property is preserved.

\medskip

\noindent\textbf{(2)} Suppose $(p,q)$ lies in a short row and a long column. Then column $p$ contains an extra cell at $(p,m+2-p)$.

For $q>1$, the label $a$ moves to $(p+1,q-1)$ and remains in a short row, and the extra cell moves to $(p+1,m+1-p)$, so the column remains long.

If $q=1$, the cell moves to $(1,m+1-p)$. Since the extra cell moves to $(p+1,m+1-p)$, the row becomes long, while the first column is short. Thus the type changes to (long row, short column).

\medskip

\noindent\textbf{(3)} Suppose $(p,q)$ lies in a long row and column $p-1$ is long. Then there is an extra cell at $(p-1,m+3-p)$.

For $q>1$, the image $(p+1,q-1)$ remains in a long row, since the extra cell structure propagates accordingly.

If $q=1$, the cell moves to $(1,m+2-p)$. The extra cell at $(p-1,m+3-p)$ moves to $(p,m+2-p)$, showing that this row is long, and the first column is also long. Hence the type becomes (long row, long column).

\medskip

\noindent\textbf{(4)} Suppose $(p,q)$ lies in a long row and column $p-1$ is short. Then there is no extra cell at $(p-1,m+3-p)$.

For $q>1$, the image $(p+1,q-1)$ remains in a long row.

If $q=1$, the cell moves to $(1,m+2-p)$. Since no extra cell moves into $(p,m+2-p)$, this row is short, while the first column is long. Hence the type becomes (short row, long column).
\end{proof}

Now we are ready to state and prove the main result of this subsection.

\begin{theorem}\label{sigma_pi^col=sigma_pi}
Let $\pi \in \mathcal{P}_n^\star$ be a periodic partition. Then
\[
  \sigma_{\pi} = \sigma_{\pi}^{\mathrm{col}} .
\]
\end{theorem}

\begin{proof}
We fix a cell $(p,q)$ of $\pi$ and follow its label under $m+1$ applications of $T_2$. By Proposition~\ref{actionT_2inverse},
the label moves one step to the right and one step down until
it reaches the bottom row, where it moves to the first column.
We use Lemma~\ref{cell movement} to determine which bottom-row
rule applies and consider the five movement types in
Remark~\ref{r1: moves of sigmapicol}.

The iteration counts in each of
Tables~\ref{Position changes for Case 1}--\ref{Position changes for Case 5} sum to $m+1$.
By Definition~\ref{actionT_2inverse=sigmapi}, the label originally
occupying the final position maps to $\pi(p,q)$ under $\sigma_\pi$.

\medskip
\noindent\textbf{Case 1. Uncolored cells:}
row $q$ and column $p$ are both short.

Both bottom-row moves use the short-row rule.
Table~\ref{Position changes for Case 1} gives the trajectory, which returns the
label to $(p,q)$. Thus
\[
\sigma_\pi(\pi(p,q))
=
\sigma_\pi^{\mathrm{col}}(\pi(p,q))
=
\pi(p,q).
\]

\begin{table}[H]
\centering
\begin{tabular}{|c|c|c|}
\hline
Initial position & Number of iterations of $T_2$ & New position \\
\hline
$(p, q)$ & $q-1$ & $(p+q-1, 1)$ \\
$(p+q-1, 1)$ & $1$ & $(1, m+2-p-q)$ \\
$(1, m+2-p-q)$ & $m+1-p-q$ & $(m+2-p-q, 1)$ \\
$(m+2-p-q, 1)$ & $1$ & $(1, p+q-1)$ \\
$(1, p+q-1)$ & $p-1$ & $(p, q)$ \\
\hline
\end{tabular}
\caption{Position changes for Case 1.}
\label{Position changes for Case 1}
\end{table}

\medskip
\noindent\textbf{Case 2. Downward move:}
row $q$ is short and column $p$ is long.

The first bottom-row move uses the short-row rule, while the
second uses the long-row rule. As shown in
Table~\ref{Position changes for Case 2}, the final position is $(p,q+1)$.
Therefore,
\[
\sigma_\pi(\pi(p,q+1))
=
\sigma_\pi^{\mathrm{col}}(\pi(p,q+1))
=
\pi(p,q).
\]

\begin{table}[H]
\centering
\begin{tabular}{|c|c|c|}
\hline
Initial position & Number of iterations of $T_2$ & New position \\
\hline
$(p, q)$ & $q-1$ & $(p+q-1, 1)$ \\
$(p+q-1, 1)$ & $1$ & $(1, m+2-p-q)$ \\
$(1, m+2-p-q)$ & $m+1-p-q$ & $(m+2-p-q, 1)$ \\
$(m+2-p-q, 1)$ & $1$ & $(1, p+q)$ \\
$(1, p+q)$ & $p-1$ & $(p, q+1)$ \\
\hline
\end{tabular}
\caption{Position changes for Case 2.}
\label{Position changes for Case 2}
\end{table}  

\medskip
\noindent\textbf{Case 3. Diagonal move:}
row $q$ and column $p-1$ are both long, with $p\geq2$.

Both bottom-row moves use the long-row rule.
Table~\ref{Position changes for Case 3} gives the final position $(p-1,q+1)$,
so
\[
\sigma_\pi(\pi(p-1,q+1))
=
\sigma_\pi^{\mathrm{col}}(\pi(p-1,q+1))
=
\pi(p,q).
\]

\begin{table}[H]
\centering
\begin{tabular}{|c|c|c|}
\hline
Initial position & Number of iterations of $T_2$ & New position \\
\hline
$(p, q)$ & $q-1$ & $(p+q-1, 1)$ \\
$(p+q-1, 1)$ & $1$ & $(1, m+3-p-q)$ \\
$(1, m+3-p-q)$ & $m+2-p-q$ & $(m+3-p-q, 1)$ \\
$(m+3-p-q, 1)$ & $1$ & $(1, p+q-1)$ \\
$(1, p+q-1)$ & $p-2$ & $(p-1, q+1)$ \\
\hline
\end{tabular}
\caption{Position changes for Case 3.}
\label{Position changes for Case 3}
\end{table}

\medskip
\noindent\textbf{Case 4. Rightward move:}
row $q$ is long and column $p-1$ is short, with $p\geq2$.

The first bottom-row move uses the long-row rule, while the
second uses the short-row rule. Table~\ref{Position changes for Case 4} therefore
gives the final position $(p-1,q)$, yielding
\[
\sigma_\pi(\pi(p-1,q))
=
\sigma_\pi^{\mathrm{col}}(\pi(p-1,q))
=
\pi(p,q).
\]

\begin{table}[H]
\centering
\begin{tabular}{|c|c|c|}
\hline
Initial position & Number of iterations of $T_2$ & New position \\
\hline
$(p, q)$ & $q-1$ & $(p+q-1, 1)$ \\
$(p+q-1, 1)$ & $1$ & $(1, m+3-p-q)$ \\
$(1, m+3-p-q)$ & $m+2-p-q$ & $(m+3-p-q, 1)$ \\
$(m+3-p-q, 1)$ & $1$ & $(1, p+q-2)$ \\
$(1, p+q-2)$ & $p-2$ & $(p-1, q)$ \\
\hline
\end{tabular}
\caption{Position changes for Case 4.}
\label{Position changes for Case 4}
\end{table} 

\medskip
\noindent\textbf{Case 5. Jump move:}
$p=1$ and $q=l_i$ or $q=h_i$.

Since row $q$ is long, the bottom-row move uses the long-row
rule. Table~\ref{Position changes for Case 5} gives the final position
$(m+2-q,1)$. Hence
\[
\sigma_\pi(\pi(m+2-q,1))
=
\sigma_\pi^{\mathrm{col}}(\pi(m+2-q,1))
=
\pi(1,q).
\]

\begin{table}[H]
\centering
\begin{tabular}{|c|c|c|}
\hline
Initial position & Number of iterations of \(T_2\) & New position \\
\hline
\((1,q)\) & \(q-1\) & \((q,1)\) \\
\((q,1)\) & \(1\) & \((1,m+2-q)\) \\
\((1,m+2-q)\) & \(m+1-q\) & \((m+2-q,1)\) \\
\hline
\end{tabular}
\caption{Position changes for Case 5.}
\label{Position changes for Case 5}
\end{table}

These cases cover every cell of $\pi$. Thus the two permutations
have the same preimage for every label, and hence
$\sigma_\pi=\sigma_\pi^{\mathrm{col}}$.
\end{proof}

\subsection{\texorpdfstring
  {The order of $\sigma_{\pi}^{col}$}{The order of the colored permutation}}\label{subsec:order of sigma}

Having established that $\sigma_{\pi}^{\mathrm{col}} = \sigma_{\pi}$, we now determine the order of this permutation. Since the order of a permutation is the least common multiple of its cycle lengths, this reduces to analyzing the cycle structure of $\sigma_{\pi}^{\mathrm{col}}$. We do this by studying the sizes of the color classes obtained with the coloring algorithm.

\begin{theorem}\label{order of sigma_pi}
Let $n \ge 1$, and let $m$ be the unique positive integer such that
\[
\Delta_m \le n < \Delta_{m+1}.
\]
Set $m' = n - \Delta_m$, and let $\pi \in \mathcal{P}_n^\star$ be periodic. Then
\[
\operatorname{ord}(\sigma_{\pi})
=
\operatorname{ord}(\sigma_{\pi}^{\mathrm{col}})
=
\begin{cases}
1, & \text{if } m' = 0, \\[4pt]
m+1, & \text{if } m' = 1, \\[4pt]
2m - m' + 3, & \text{if } 2 \le m' \le m .
\end{cases}
\]
\end{theorem}

\begin{proof}
Let $m' = n - \Delta_m$. We consider three cases.

\medskip
\noindent\textbf{Case 1: $m' = 0$.}
Then $\pi$ is the full triangular partition of size $\Delta_m$. In this case, $\sigma_{\pi}^{\mathrm{col}}$, and hence also $\sigma_\pi$, is the identity permutation. Therefore,
\[
\operatorname{ord}(\sigma_\pi) = 1.
\]

\medskip
\noindent\textbf{Case 2: $m' = 1$.}
There is exactly one extra cell, say at location $(i,h)$ in the $(m+1)$-st diagonal. The coloring algorithm colors this cell together with $i-1$ cells to its left and $h-1$ cells below. The total number of colored cells is
\[
1 + (i-1) + (h-1) = i + h - 1 = m + 1.
\]
These labels form a single cycle, so
$\operatorname{ord}(\sigma_\pi)=m+1$.

\medskip
\noindent\textbf{Case 3: \(2 \le m' \le m\).}
Let the extra cells on $\diag[m+1]$ be
\[
E=\bigl((i_1,h_1),\dots,(i_{m'},h_{m'})\bigr),
\qquad i_1<\cdots<i_{m'},
\]
so $i_k+h_k=m+2$ and hence $h_1>\cdots>h_{m'}$.
For $0\le t\le m'-1$ set
\[
X^t(E)=\bigl((i_1,h_{1+t}),\dots,(i_{m'-t},h_{m'})\bigr).
\]

Fix $1\le r\le\lceil m'/2\rceil$.  The color $C_r$ consists of the two components
$X^{r-1}(E)$ and $X^{m'-r}(E)$ (coinciding if $m'$ is odd and
$r=\tfrac{m'+1}{2}$).
Write
\[
X^t(E)=((c_1,d_1),\dots,(c_{m'-t},d_{m'-t})).
\]

Leftward propagation in row \(d_j\) colors \(c_j - c_{j-1} - 1\) cells (with \(c_0 = 0\)), while downward propagation in column \(c_j\) colors \(d_j - d_{j+1} - 1\) cells (with \(d_{m'-t+1} = 0\)). Thus, the total number of cells colored by \(C_r\) in this component, starting from the cells in \(X^t(E)\), is

\[
 \sum_{j=1}^{m'-t} \bigl(1 + c_j - c_{j-1} - 1 + d_j - d_{j+1} - 1\bigr) 
= c_{m'-t} + d_1 - (m'-t)
= i_{m'-t} + h_{1+t} - m' + t.
\]

Suppose first that the two chains of $C_r$ are distinct. Applying the preceding calculation to the two chains of $C_r$,
\begin{align*}
|C_r|
&=\bigl(i_{m'-r+1}+h_r-m'+r-1\bigr)
  +\bigl(i_r+h_{m'-r+1}-r\bigr)\\
&=\bigl(i_{m'-r+1}+m+2-i_r-m'+r-1\bigr)
  +\bigl(i_r+m+2-i_{m'-r+1}-r\bigr)\\
&=2m-m'+3 .
\end{align*}
Thus every color class consisting of two distinct chains has
$2m-m'+3$ cells. If $m'$ is odd and
$r=(m'+1)/2$, the two chains coincide, so the final color class
has $(2m-m'+3)/2$ cells. If $m'$ is even, every color class has $2m-m'+3$ cells.  Since $m'\geq2$, at least one cycle has length $2m-m'+3$, and every other cycle length divides this number.
Therefore,
\[
\operatorname{ord}(\sigma_\pi^{\mathrm{col}})
=2m-m'+3.
\]
The equality $\sigma_\pi=\sigma_\pi^{\mathrm{col}}$ now gives
the stated formula.
\end{proof}

\begin{corollary}\label{cor:orbit-size-T2}
    Let $\pi \in \mathcal{P}_n^\star$ be a periodic point with $n = \Delta_m + m'$ and $0 \le m' \le m$. Then the orbit size of $\pi$ under $T_2$ divides
    \[
\begin{cases}
m+1, & \text{if } m' = 0, \\[2mm]
(m+1)^2, & \text{if } m' = 1, \\[2mm]
(m+1)(2m - m' + 3), & \text{if } 2 \le m' \le m.
\end{cases}
\]

   Moreover, if the period of the underlying shape \(S(\pi)\) under \(T\) is \(\ell = m+1\), then ``divides'' may be replaced by ``equals''.
\end{corollary}

\begin{proof}
The claim follows from substituting the values of $\operatorname{ord}(\sigma_\pi)$ from
Theorem~\ref{order of sigma_pi} into Theorem~\ref{thm:maximal orbit size}.
\end{proof}

\medskip

The following examples show that the bound may or may not be attained when the underlying shape has period less than $m+1$.

\begin{example}[$n=18$, $m=5$, $m'=3$]
\label{ex:all-four-orbits-18}
There are four periodic orbits in $\mathcal{P}_{18}$ under $T$,
corresponding to necklaces of length $6$ with $3$ black beads. Figure~\ref{fig:four-orbits-18} shows a labeled representative of each underlying orbit, together with its shape period $\ell$.

\begin{figure}[H]
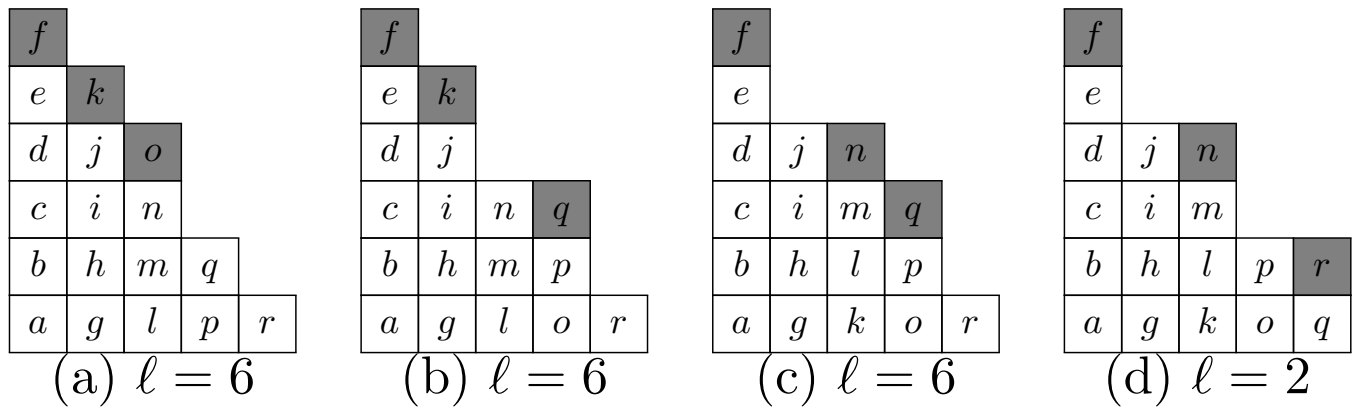

\centering
\renewcommand{\arraystretch}{0.8}
\begin{tabular}{cccc}
\scalebox{0.7}{\begin{ytableau}
*(gray)f\empty & \none & \none & \none & \none \\
e\empty & *(gray)k\empty & \none & \none & \none \\
d\empty & j\empty & *(gray)o\empty & \none & \none \\
c\empty & i\empty & n\empty & \none & \none \\
b\empty & h\empty & m\empty & q\empty & \none \\
a\empty & g\empty & l\empty & p\empty & r\empty
\end{ytableau}}
&
\scalebox{0.7}{\begin{ytableau}
*(gray)f\empty & \none & \none & \none & \none \\
e\empty & *(gray)k\empty & \none & \none & \none \\
d\empty & j\empty & \none & \none & \none \\
c\empty & i\empty & n\empty & *(gray)q\empty & \none \\
b\empty & h\empty & m\empty & p\empty & \none \\
a\empty & g\empty & l\empty & o\empty & r\empty
\end{ytableau}}
&
\scalebox{0.7}{\begin{ytableau}
*(gray)f\empty & \none & \none & \none & \none \\
e\empty & \none & \none & \none & \none \\
d\empty & j\empty & *(gray)n\empty & \none & \none \\
c\empty & i\empty & m\empty & *(gray)q\empty & \none \\
b\empty & h\empty & l\empty & p\empty & \none \\
a\empty & g\empty & k\empty & o\empty & r\empty
\end{ytableau}}
&
\scalebox{0.7}{\begin{ytableau}
*(gray)f\empty & \none & \none & \none & \none \\
e\empty & \none & \none & \none & \none \\
d\empty & j\empty & *(gray)n\empty & \none & \none \\
c\empty & i\empty & m\empty & \none & \none \\
b\empty & h\empty & l\empty & p\empty & *(gray) r\empty \\
a\empty & g\empty & k\empty & o\empty & q\empty
\end{ytableau}}
\\[4pt]
(a) $\ell=6$ & (b) $\ell=6$ & (c) $\ell=6$  & (d) $\ell=2$  
\end{tabular}
\caption{Labeled configurations whose underlying shapes represent the four distinct periodic orbits in $\mathcal{P}_{18}$ under $T$.
Extra cells (the $m'=3$ cells on $\operatorname{diag}[6]$) are indicated in gray. The underlying shapes in (a)--(c) have period $\ell=6$ under $T$, while the shape in (d) has period $\ell=2$.}
\label{fig:four-orbits-18}
\end{figure}

By Corollary~\ref{cor:orbit-size-T2}, since $\ell = m+1$, the first three orbits under the action of $T_2$ have size $(m+1)(2m-m'+3)=60$, while the fourth orbit has size dividing 60.
However, direct computation shows that the fourth orbit also has size 60.
\end{example}
\medskip

\begin{example}[$n=8$, $m=3$, $m'=2$]
\label{ex:m-mprime-1}
There are two periodic orbits in $\mathcal{P}_{8}$ under $T$,
corresponding to necklaces of length $4$ with $2$ black beads. Figure~\ref{fig:two-orbits-8} shows a labeled representative of each underlying orbit, together with its shape period $\ell$.

\begin{figure}[H]
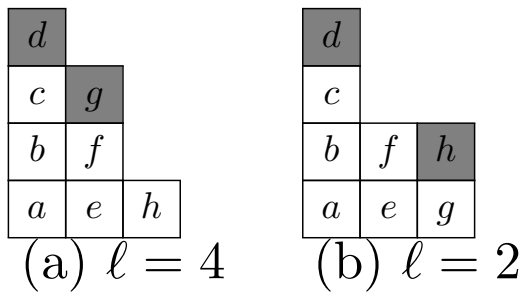

\centering
\renewcommand{\arraystretch}{0.8}
\begin{tabular}{cc}
\scalebox{0.7}{\begin{ytableau}
*(gray) d & \none & \none & \none \\
c & *(gray) g & \none & \none \\
b & f & \none & \none \\
a & e & h & \none
\end{ytableau}}
&
\scalebox{0.7}{\begin{ytableau}
*(gray) d & \none & \none & \none \\
c & \none & \none & \none \\
b & f & *(gray) h & \none \\
a & e & g & \none
\end{ytableau}}
\\[4pt]
(a) $\ell=4$ & (b) $\ell=2$  
\end{tabular}
\caption{Labeled configurations whose underlying shapes represent the two distinct periodic orbits in $\mathcal{P}_8$ under $T$.
Extra cells on $\operatorname{diag}[4]$ are shaded gray and labeled $d,g$ in (a) and $d,h$ in (b). The underlying shapes in (a) and (b)
have periods $\ell=4$ and $\ell=2$ under $T$, respectively.}
\label{fig:two-orbits-8}
\end{figure}

Here $\operatorname{ord}(\sigma_\pi)=2m-m'+3=7$. By Corollary~\ref{cor:orbit-size-T2}, since $\ell = m+1$ for orbit (a), its size under $T_2$ equals $(m+1)\operatorname{ord}(\sigma_\pi)=28$; for orbit (b) with $\ell=2$, the orbit size divides $28$, and direct computation shows it is $14$.
\end{example}

\noindent
These examples illustrate the bound of Corollary~\ref{cor:orbit-size-T2}: the orbit size under $T_2$ always divides $(m+1)\operatorname{ord}(\sigma_\pi)$, and sometimes attains it.  
In the next section, we determine the exact orbit sizes and prove that the upper bound $(m+1)\operatorname{ord}(\sigma_\pi)$ from Corollary~\ref{cor:orbit-size-T2} is achieved in most cases, even when $\ell \neq m+1$, and when it is not, the orbit size is exactly half of that bound.

\section{\texorpdfstring
  {Exact sizes of orbits under $T_2$}{Exact sizes of orbits under T2}}\label{subsec:orbit-size}

In this section, we determine the exact orbit sizes of periodic partitions under $T_2$ and characterize when the bound $(m+1)\operatorname{ord}(\sigma_\pi)$ from
Corollary~\ref{cor:orbit-size-T2} is attained. As the preceding examples show, equality can hold even when the period $\ell$ of the underlying shape is less than $m+1$. We prove that, for $m\geq3$, the bound is attained except when $m'=m-1$ and $\ell=(m+1)/2$, in which case the orbit size is exactly half of the bound. To establish these results, we analyze the permutation induced by one full period of the underlying shape, treating the colored and uncolored labels separately.

We begin by introducing notation for the colored and uncolored labels.

\begin{definition}\label{def:label-sets} Let $\pi\in\mathcal{P}_n^\star$ be a periodic partition whose cells are colored using the algorithm in Section~\ref{subsec:algo}. Define:

\begin{itemize}
    \item $\mathcal{L}_{\pi,C}$ as the set of labels of colored cells of $\pi$,
    \item $\mathcal{L}_{\pi,U}$ as the set of labels of uncolored cells of $\pi$,
    \item for each color $C_i$ ($i=1,\dots,\lceil m'/2 \rceil$), $\mathcal{L}_{\pi,C_i}$ as the set of labels of cells colored with $C_i$.
\end{itemize}
\end{definition}
Note that
\[
\mathcal{L}_\pi = \mathcal{L}_{\pi,C} \sqcup \mathcal{L}_{\pi,U}
\quad \text{and} \quad
\mathcal{L}_{\pi,C} = \bigsqcup_{i=1}^{\lceil m'/2 \rceil} \mathcal{L}_{\pi,C_i}.
\]

We now introduce the permutation of all labels induced by one full period of the underlying shape.

\begin{definition}\label{def:tau_pi}
Let \(\pi \in \mathcal{P}_n^\star\) be a periodic partition. The permutation \(\tau_\pi: \mathcal{L}_\pi \rightarrow \mathcal{L}_\pi\) is defined by
\[
\tau_\pi(\pi(i,j)) = (T_2^{\,\ell}\pi)(i,j),
\]
where $\ell$ is the period of the underlying unlabeled partition under \(T\).
\end{definition}

\begin{lemma}\label{lem:tau-cycle-structure-combined}
Let \(\pi \in \mathcal{P}_n^\star\) be a periodic point with \(n = \Delta_m + m'\), \(0 \leq m' \leq m\).  
Let \(\ell\) be the period of the underlying unlabeled partition under \(T\) and set \(r = (m+1)/\ell\).  
For each color class \(C_i\) (\(i = 1,\dots, \lceil m'/2 \rceil\)) the following hold:

\begin{enumerate}[label=(\alph*)]
    \item The set \(\mathcal{L}_{\pi,C_i}\) is invariant under \(\tau_\pi\), i.e.,
    \[
    \tau_\pi(\mathcal{L}_{\pi,C_i}) = \mathcal{L}_{\pi,C_i}.
    \]

    \item The restriction \(\tau_{\pi,C_i} := \tau_\pi|_{\mathcal{L}_{\pi,C_i}}\) is a permutation of \(\mathcal{L}_{\pi,C_i}\) and satisfies
    \[
    \tau_{\pi,C_i}^{\,r} = \sigma_{\pi,C_i},
    \]
where $\sigma_{\pi,C_i}:=\left.\sigma_\pi\right|_{\mathcal{L}_{\pi,C_i}}$.
    \item The permutation \(\tau_{\pi,C_i}\) is a single cycle on \(\mathcal{L}_{\pi,C_i}\).
\end{enumerate}
\end{lemma}

\begin{proof}

\begin{enumerate}[label=(\alph*)]
    \item We first show that a label retains its color after one
application of $T_2$, when the resulting configuration is colored
by the algorithm.

Let $(p,q)$ be a colored cell, and let
\[
d=\#\{s:p\leq i_s\leq m+2-q\},
\]
where $i_1,\ldots,i_{m'}$ are the column indices of the extra cells.
The color of $(p,q)$ is $C_{\min(d,m'+1-d)}$. Indeed, at a cell
of $X^{d-1}(E)$ the indicated interval contains exactly $d$
extra-cell column indices, and propagation along the corresponding
chain does not change this number.

If $q>1$, Proposition~\ref{actionT_2inverse} sends $(p,q)$ to $(p+1,q-1)$.
The interval and the extra-cell positions shift together, so $d$
is unchanged.

Suppose $q=1$. If the bottom row is long, the label moves to
$(1,m+2-p)$, and the number of extra-cell column indices in the
new interval is $m'+1-d$. If the bottom row is short, then column
$p$ must be long because the cell is colored. The label moves
to $(1,m+1-p)$, and the new count is again $m'+1-d$.
In both cases the resulting cell is colored, and
\[
\min(d,m'+1-d)=\min(m'+1-d,d).
\]
Thus the label retains its color.

Iterating this observation $\ell$ times, and using the fact that
$T_2^\ell(\pi)$ has the same shape as $\pi$, shows that
$\tau_\pi(\mathcal{L}_{\pi,C_i})=\mathcal{L}_{\pi,C_i}$.

\item Therefore, the restriction
\[
\tau_{\pi,C_i}:=\tau_\pi\big|_{\mathcal{L}_{\pi,C_i}}
\]
is a permutation of \(\mathcal{L}_{\pi,C_i}\). Since \(m+1=\ell r\), by Definitions  \ref{actionT_2inverse=sigmapi} and \ref{def:tau_pi},
\[
\tau_\pi^{\,r}=\sigma_\pi.
\]
By Theorem~\ref{sigma_pi^col=sigma_pi}, we also have
\(\sigma_\pi=\sigma_\pi^{\mathrm{col}}\). Restricting this identity to
\(\mathcal{L}_{\pi,C_i}\) gives
\[
\tau_{\pi,C_i}^{\,r}
=\sigma_\pi\big|_{\mathcal{L}_{\pi,C_i}}
=\sigma_\pi^{\mathrm{col}}\big|_{\mathcal{L}_{\pi,C_i}}
=\sigma_{\pi,C_i}.
\]

\item By construction, each \(\sigma_{\pi,C_i}\) is a single cycle. Since \(\sigma_{\pi,C_i}=\tau_{\pi,C_i}^{\,r}\), a power of \(\tau_{\pi,C_i}\) is a single cycle. Hence \(\tau_{\pi,C_i}\) itself must be a single cycle on \(\mathcal{L}_{\pi,C_i}\).

This completes the proof.
\end{enumerate}
\end{proof}

To combine the contributions of the colored and uncolored labels to the orbit size, we will also need the following coprimality condition.

\begin{lemma}\label{lem:gcd-n-r}
Let $\pi \in \mathcal{P}_n^\star$ be a periodic partition with $n = \Delta_m + m'$ and $0 \leq m' \leq m$. Let $\ell$ denote the period of the underlying unlabeled partition under $T$, so that $\ell \mid (m+1)$, and set $r = (m+1)/\ell$. Then
\[
\gcd\bigl(\operatorname{ord}(\sigma_\pi),\, r\bigr) = 1.
\]
\end{lemma}

\begin{proof}
We consider three cases based on the value of $m'$.

\medskip

\noindent\textbf{Case \(m' = 0\).}
Then, by Theorem~\ref{order of sigma_pi}, \(\operatorname{ord}(\sigma_\pi) = 1\), so \(\gcd(\operatorname{ord}(\sigma_\pi),r) = 1\).

\medskip

\noindent\textbf{Case \(m' = 1\).}
Then, by Theorem~\ref{order of sigma_pi}, \(\operatorname{ord}(\sigma_\pi) = m+1\). Since there is a single extra cell, the underlying unlabeled partition has full period \(\ell = m+1\), hence \(r = 1\). Thus \(\gcd(\operatorname{ord}(\sigma_\pi),r) = 1\).

\medskip

\noindent\textbf{Case \(2 \le m' \le m\).}
Then, by Theorem~\ref{order of sigma_pi}, we have \(\operatorname{ord}(\sigma_\pi) = 2m - m' + 3\). The binary word recording the occupied positions on $\operatorname{diag}[m+1]$ consists of $r$ repetitions of a block
of length $\ell$. Therefore its number $m'$ of occupied positions
is divisible by $r$. Write \(m' = ra\) for some integer \(a \ge 1\). Since \(m = \ell r - 1\), we compute
\[
\operatorname{ord}(\sigma_\pi) = 2m - m' + 3 = 2(\ell r - 1) - ra + 3 =  r(2\ell - a) + 1.
\]
Hence \(\operatorname{ord}(\sigma_\pi) \equiv 1 \pmod{r}\), which implies \(\gcd(\operatorname{ord}(\sigma_\pi),r) = 1\).
\end{proof}

We are now ready to state the main theorem that gives the exact orbit sizes under $T_2$.

\begin{theorem}\label{thm:orbit-size}
Let \(\pi \in \mathcal{P}_n^\star\) be a periodic point under the action \(T_2\), with \(n \ge 1\) and \(n = \Delta_m + m'\), \(0 \le m' \le m\). Let $\ell$ be the period of the underlying unlabeled partition under \(T\). Then, the orbit size of \(\pi\) under \(T_2\) is
\[
|O_{T_2}(\pi)| =  \begin{cases}
(m+1)\cdot \operatorname{ord}(\sigma_\pi), & \text{if } m-m' \ge 2, \\[2mm]
\ell \cdot \operatorname{ord}(\sigma_\pi), & \text{otherwise}.
\end{cases}
\]

\end{theorem}

\begin{proof}
Let $L=\operatorname{ord}(\sigma_\pi)$ and $r=(m+1)/\ell$, and
let $\tau_\pi$ be the permutation from Definition~\ref{def:tau_pi}. Since
$m+1=\ell r$, Definitions~\ref{actionT_2inverse=sigmapi} and~\ref{def:tau_pi} give
$\tau_\pi^r=\sigma_\pi$.

By Lemma~\ref{lem:tau-cycle-structure-combined}, $\tau_\pi$ preserves each colored label set
$\mathcal{L}_{\pi,C_i}$. Consequently, it preserves their union
$\mathcal{L}_{\pi,C}$ and its complement $\mathcal{L}_{\pi,U}$.
Define the restrictions
\[
\tau_c=\left.\tau_\pi\right|_{\mathcal{L}_{\pi,C}},
\qquad
\tau_u=\left.\tau_\pi\right|_{\mathcal{L}_{\pi,U}},
\]
and let $m_c=\operatorname{ord}(\tau_c)$ and
$m_u=\operatorname{ord}(\tau_u)$. Since these restrictions act
on disjoint label sets,
\[
\operatorname{ord}(\tau_\pi)=\operatorname{lcm}(m_c,m_u).
\]
Moreover, $\sigma_\pi$ fixes every uncolored label, so restricting
$\tau_\pi^r=\sigma_\pi$ to $\mathcal{L}_{\pi,U}$ gives
$\tau_u^r=\mathrm{id}$. Hence $m_u\mid r$.

Every return time of $\pi$ under $T_2$ is a multiple of $\ell$,
since its underlying shape has period $\ell$. Because the shape
repeats every $\ell$ steps and the rules of $T_2$ do not depend
on the labels, $T_2^{k\ell}(\pi)=\pi$ if and only if
$\tau_\pi^k=\mathrm{id}$. Therefore,
\[
|O_{T_2}(\pi)|
=\ell\,\operatorname{ord}(\tau_\pi)
=\ell\,\operatorname{lcm}(m_c,m_u).
\]

By Lemma~\ref{lem:tau-cycle-structure-combined}, each restriction $\tau_{\pi,C_i}$ is a single cycle on the same label set as $\sigma_{\pi,C_i}$, so these cycles have
the same length. Since $\sigma_\pi$ fixes all uncolored labels, taking the least common multiple of the colored cycle lengths gives
$m_c=\operatorname{ord}(\tau_c)=\operatorname{ord}(\sigma_\pi)=L$.
This also holds when $m'=0$, since both permutations are identities.

Next, we determine $m_u = \operatorname{ord}(\tau_u)$. Suppose first that there is no more than one uncolored cell. Then $\tau_u$ must be the identity on at most one label, hence $m_u = 1$.

Suppose now that there are at least three uncolored cells.

Suppose now that there at least two uncolored cells. Let $k=m-m'$. Let $A\subseteq\{1,\ldots,m+1\}$ be the set of column indices
of the unoccupied positions on $\operatorname{diag}[m+1]$.
Then $|A|=m+1-m'=k+1\geq3$.

The uncolored labels are in bijection with the two-element
subsets of $A$ via
\[
\pi(p,q)\longleftrightarrow\{p,m+2-q\}.
\]
Indeed, an uncolored cell lies in a short column and a short row,
and $p<m+2-q$. Conversely, for $a<b$ in $A$, the cell
$(a,m+2-b)$ is uncolored.

By Proposition~\ref{actionT_2inverse}, if $q>1$, the cell $(p,q)$ moves to
$(p+1,q-1)$, so its associated pair becomes
$\{p+1,m+3-q\}$. If $q=1$, the bottom row must be short because
the cell is uncolored. Hence the second branch of Proposition~\ref{actionT_2inverse}
applies, sending $(p,1)$ to $(1,m+1-p)$, whose associated pair
is $\{1,p+1\}$. In both cases, the associated unordered pair is
obtained by adding $1$ modulo $m+1$ to each of its indices,
with residues represented by $1,\ldots,m+1$.  Since the shape returns after
$\ell$ steps, addition of $\ell$ modulo $m+1$ permutes $A$.
Its induced action on the two-element subsets of $A$ is the
inverse of $\tau_u$.

Suppose $\tau_u^d=\mathrm{id}$. Then addition of $\ell d$ modulo
$m+1$ fixes every two-element subset of $A$. For any $a\in A$,
choose distinct $b,c\in A\setminus\{a\}$. Since this translation
fixes both $\{a,b\}$ and $\{a,c\}$, it fixes their intersection
$\{a\}$. Thus $\ell d\equiv0\pmod{m+1}$, which implies $r\mid d$. Together with $\tau_u^r=\mathrm{id}$, this proves that
$m_u=\operatorname{ord}(\tau_u)=r$.

We have shown that $m_c=L$, while $m_u=r$ if $m-m'\geq2$
and $m_u=1$ otherwise. By Lemma~\ref{lem:gcd-n-r}, $\gcd(L,r)=1$.
Consequently,
\[
|O_{T_2}(\pi)|
=\ell\,\operatorname{lcm}(m_c,m_u)
=
\begin{cases}
\ell Lr=(m+1)L, & \text{if }m-m'\geq2,\\
\ell L, & \text{otherwise}.
\end{cases}
\]
Since $L=\operatorname{ord}(\sigma_\pi)$, this proves the theorem.
\end{proof}

Combining Theorem~\ref{thm:orbit-size} with the explicit formula for
\(\operatorname{ord}(\sigma_\pi)\) from Theorem~\ref{order of sigma_pi},
we obtain the following closed-form expressions for orbit sizes when \(m\ge 3\).

\begin{theorem}

Let \(\pi \in \mathcal{P}_n^\star\) be a periodic point under the action \(T_2\),
where \(n=\Delta_m+m'\) with \(m\ge 3\) and \(0\le m'\le m\).
Let \(\ell\) be the period of the underlying unlabeled partition under \(T\). Then, the orbit size of \(\pi\) under \(T_2\) is
\[
\bigl|O_{T_2}(\pi)\bigr|=
\begin{cases}
m+1, 
& \text{if } m'=0 \ (n=\Delta_m), \\[2mm]

(m+1)^2, 
& \text{if } m'=1 \ (n=\Delta_m+1), \\[2mm]

(m+1)(2m-m'+3), 
& \text{if } 2\le m'\le m-2, \\[2mm]

(m+1)(m+4)\ \text{or}\ \dfrac{(m+1)(m+4)}{2}^{\ast},
& \text{if } m'=m-1 \ (n=\Delta_m+m-1), \\[3mm]

(m+1)(m+3), 
& \text{if } m'=m \ (n=\Delta_m+m).
\end{cases}
\]
\medskip

\noindent
$\ast$~The second value occurs exactly when $\ell=(m+1)/2$, equivalently when the two white beads in the associated necklace are opposite each other. This is possible only when $m$ is odd.
\end{theorem}

\begin{proof}
For $0\leq m'\leq m-2$, we have $m-m'\geq2$, so the stated
orbit sizes follow directly from Theorems~\ref{order of sigma_pi} and~\ref{thm:orbit-size}.

Suppose $m'=m-1$. These theorems give
$|O_{T_2}(\pi)|=\ell(m+4)$. The associated necklace has length
$m+1$ and exactly two white beads. Any nonidentity rotation
preserving its bead coloring must interchange the two white
beads and is therefore a half-turn. Hence $\ell=m+1$ unless
the white beads are opposite each other, in which case
$\ell=(m+1)/2$. The latter is possible only when $m$ is odd,
giving the two stated orbit sizes.

Finally, if $m'=m$, the associated necklace has exactly one
white bead and therefore has full period $\ell=m+1$.
Since $\operatorname{ord}(\sigma_\pi)=m+3$, Theorem~5.5 gives
$|O_{T_2}(\pi)|=(m+1)(m+3)$.
\end{proof}

\section{Conjugation symmetries}\label{sec:conjugation}

In Section~\ref{sec:actions}, we established how the reversal operation interacts with our four actions $T_i$, $i=1,2,3,4$: it conjugates $T_1$ to $T_4$ and $T_2$ to $T_3$. Conjugation
of Ferrers diagrams is another standard operation on partitions with distinguishable cells, so it is natural to ask how it
interacts with these actions. Here we establish such relations for $T_1$ and $T_2$ on periodic
partitions and use reversal to obtain corresponding relations for $T_3$ and $T_4$.

\begin{definition}\label{def:conjugation}
Let \(\pi \in \mathcal{P}_n^\star\). The \emph{conjugate} of \(\pi\), denoted \(\conj(\pi)\), is the partition obtained by reflecting the Ferrers diagram of \(\pi\) across its main diagonal. Equivalently,
\[
\operatorname{conj}(\pi)(i,j)=\pi(j,i)
\]
for every cell $(i,j)$ in the Ferrers diagram of
$\operatorname{conj}(\pi)$. Thus, the label occupying row \(i\) and column \(j\) in \(\pi\) occupies row \(j\) and column \(i\) in \(\conj(\pi)\).
Note that conjugation is an involution; that is, \(\conj(\conj(\pi)) = \pi\) for every \(\pi\in\mathcal{P}_n^\star\).
\end{definition}
Since conjugation preserves each diagonal $\operatorname{diag}[k]$ setwise, Theorem~\ref{periodic points under T} and Proposition~\ref{prop:periodic shape} imply that it preserves
periodicity. We examine how the actions $T_1$ and $T_2$ interact with conjugation.

\begin{proposition}\label{prop:conj-relation}
For every periodic \(\pi\in\mathcal{P}_n^\star\), the following diagrams commute:
\[
\begin{tikzcd}
\pi \ar[r,"T_1"] \ar[d,"\conj"'] & T_1(\pi) \ar[d,"\conj"] \\
\conj(\pi) & \conj(T_1(\pi)) \ar[l,"T_1"']
\end{tikzcd}
\qquad\qquad
\begin{tikzcd}
\pi \ar[r,"T_2"] \ar[d,"\conj"'] & T_2(\pi) \ar[d,"\conj"] \\
\conj(\pi) & \conj(T_2(\pi)) \ar[l,"T_2"']
\end{tikzcd}
\]
Equivalently,
\[
T_1\circ\conj\circ T_1=\conj,
\qquad
T_2\circ\conj\circ T_2=\conj.
\]
\end{proposition}

\begin{proof}
Let \(\pi \in \mathcal{P}_n^\star\) be periodic. Under conjugation, the label \(x=\pi(i,j)\) moves  to position \((j,i)\) in \(\conj(\pi)\).

\smallskip
\noindent\textbf{The action \(T_1\).} Under $T_1$, the label $x$ moves to $(1,i)$ if $j=1$, and to $(i+1,j-1)$ if $j>1$. After conjugation, these positions become
\[
(i,1)\qquad\text{if }j=1,
\qquad\text{and}\qquad
(j-1,i+1)\qquad\text{if }j>1.
\]
Applying \(T_1\) once more, the entry at \((i,1)\) moves to \((1,i)\), whereas the cell at \((j-1,i+1)\) moves to \((j,i)\). Thus in either case the final position is \((j,i)\), exactly the position of \(x\) in \(\conj(\pi)\). Therefore
\[
T_1\circ\conj\circ T_1=\conj.
\]

\smallskip
\noindent\textbf{The action \(T_2\).}
Let \(k\) denote the number of columns of \(\pi\), equivalently the length of its bottom row. Under \(T_2\), an entry with \(j>1\) again moves to \((i+1,j-1)\), while a bottom row entry at \((i,1)\) is placed in the new first column at position \((1,k+1-i)\). After conjugation, these become
\[
(k+1-i,1)\qquad\text{if }j=1,
\qquad\text{and}\qquad
(j-1,i+1)\qquad\text{if }j>1.
\]
Since $\pi$ is periodic, the newly formed column of $T_2(\pi)$ has length $k$ and is at least as long as every remaining column.
Thus $T_2(\pi)$ has exactly $k$ rows, so the bottom row of $\operatorname{conj}(T_2(\pi))$ has length $k$. Hence applying \(T_2\) again sends the entry at \((k+1-i,1)\) to
\(
(1,k+1-(k+1-i))=(1,i),
\)
while the entry at \((j-1,i+1)\) moves to \((j,i)\). Thus the final position is again \((j,i)\), which is precisely the position of \(x\) in \(\conj(\pi)\). Therefore
\[
T_2\circ\conj\circ T_2=\conj.
\]

Since this holds for every entry \(x\) of \(\pi\), the result follows.
\end{proof}

On the set of periodic partitions, Proposition~\ref{prop:conj-relation} can equivalently
be written as
\[
\operatorname{conj}\circ T_i
=T_i^{-1}\circ\operatorname{conj},
\qquad i=1,2.
\]
Thus conjugation sends each periodic orbit to a periodic orbit
of the same size and reverses the direction of the action.

\begin{remark}
Combining Propositions~\ref{prop:rev-relation} and~\ref{prop:conj-relation} gives corresponding identities
for $T_3$ and $T_4$. On the set of periodic partitions,
\[
T_i\circ
(\operatorname{rev}\circ\operatorname{conj}\circ\operatorname{rev})
\circ T_i
=
\operatorname{rev}\circ\operatorname{conj}\circ\operatorname{rev},
\qquad i=3,4.
\]
Indeed, these identities follow by composing the identities in
Proposition~\ref{prop:conj-relation} on both sides with $\operatorname{rev}$ and using
Proposition~\ref{prop:rev-relation}. Thus, the involution
$\operatorname{rev}\circ\operatorname{conj}\circ\operatorname{rev}$
reverses the direction of the action on periodic orbits of
$T_3$ and $T_4$, just as $\operatorname{conj}$ does for
$T_1$ and $T_2$.
\end{remark}

\begin{remark}
Proposition~\ref{prop:conj-relation} does not hold for arbitrary
(non‑periodic) partitions.  
For example, let \(n=3\) and consider the one‑column partition  
\(\pi = ((a,b,c)) \in \mathcal{P}_3^\star\) (where the cells are listed from
bottom to top).  
Then \(\conj(\pi) = ((a),(b),(c))\); in particular \(\conj(\pi)\) has three
columns of size \(1\).  
On the other hand,
\[
T_2(\pi) = ((b,c))\cup (a) = ((b,c),(a)),
\qquad
\conj\bigl(T_2(\pi)\bigr) = ((b,a),(c)),
\]
and finally
\[
T_2\bigl(\conj(T_2(\pi))\bigr) = ((c,b),(a)) \neq \conj(\pi).
\]

Thus the identity for $T_2$ fails for this nonperiodic partition.
Similarly,
\[
T_1\bigl(\operatorname{conj}(T_1(\pi))\bigr)
=((b,c),(a))\neq\operatorname{conj}(\pi).
\]
Consequently, the periodicity assumption in Proposition~\ref{prop:conj-relation} cannot be omitted.
\end{remark}

\bibliographystyle{plain}
\bibliography{refs}

@article {hopkins2024henrik,
    AUTHOR = {Hopkins, Brian},
     TITLE = {Henrik {E}riksson's {B}ulgarian solitaire variant},
   JOURNAL = {Integers},
  FJOURNAL = {Integers. Electronic Journal of Combinatorial Number Theory},
    VOLUME = {24A},
      YEAR = {2024},
     PAGES = {Paper No. A9, 14},
      ISSN = {1553-1732},
   MRCLASS = {05C57 (11P81)},
  MRNUMBER = {4750801},
MRREVIEWER = {Warren\ P.\ Johnson},
       DOI = {10.5281/zenodo.11352790},
       URL = {https://doi.org/10.5281/zenodo.11352790},
}

@article {brandt1982cycles,
    AUTHOR = {Brandt, J{\o}rgen},
     TITLE = {Cycles of partitions},
   JOURNAL = {Proc. Amer. Math. Soc.},
  FJOURNAL = {Proceedings of the American Mathematical Society},
    VOLUME = {85},
      YEAR = {1982},
    NUMBER = {3},
     PAGES = {483--486},
      ISSN = {0002-9939,1088-6826},
   MRCLASS = {05A17 (10A45)},
  MRNUMBER = {656129},
MRREVIEWER = {Bruce\ Sagan},
       DOI = {10.2307/2043873},
       URL = {https://doi-org.libproxy.clemson.edu/10.2307/2043873},
}

@incollection{gardner1997bulgarian,
  title={Bulgarian solitaire and other seemingly endless tasks},
  author={Gardner, Martin},
  booktitle={The Last Recreations: Hydras, Eggs, and Other Mathematical Mystifications},
  pages={27--43},
  year={1997},
  publisher={Springer}
}

@article{toom1981,
  author  = {Toom, A. L.},
  title   = {Problem Solution {M655}},
  journal = {Kvant},
  volume  = {12},
  number  = {7},
  year    = {1981},
  pages   = {28--30},
  }

@article {hopkins2012,
    AUTHOR = {Hopkins, Brian},
     TITLE = {30 years of {B}ulgarian solitaire},
   JOURNAL = {College Math. J.},
  FJOURNAL = {The College Mathematics Journal},
    VOLUME = {43},
      YEAR = {2012},
    NUMBER = {2},
     PAGES = {135--140},
      ISSN = {0746-8342,1931-1346},
   MRCLASS = {11P81 (91A05)},
  MRNUMBER = {2897476},
       DOI = {10.4169/college.math.j.43.2.135},
       URL = {https://doi-org.libproxy.clemson.edu/10.4169/college.math.j.43.2.135},
}

@article {hart2012,
    AUTHOR = {Hart, Therese A. and Khan, Gabriel and Khan, Mizan R.},
     TITLE = {Revisiting {T}oom's proof of {B}ulgarian solitaire},
   JOURNAL = {Ann. Sci. Math. Qu\'{e}bec},
  FJOURNAL = {Annales des Sciences Math\'{e}matiques du Qu\'{e}bec},
    VOLUME = {36},
      YEAR = {2012},
    NUMBER = {2},
     PAGES = {477--486},
      ISSN = {0224-7992},
   MRCLASS = {05A17 (05A15)},
  MRNUMBER = {3115424},
MRREVIEWER = {Christian Krattenthaler},
}

@article {igusa1985,
    AUTHOR = {Igusa, Kiyoshi},
     TITLE = {Solution of the {B}ulgarian solitaire conjecture},
   JOURNAL = {Math. Mag.},
  FJOURNAL = {Mathematics Magazine},
    VOLUME = {58},
      YEAR = {1985},
    NUMBER = {5},
     PAGES = {259--271},
      ISSN = {0025-570X,1930-0980},
   MRCLASS = {00A08 (90D99)},
  MRNUMBER = {810147},
MRREVIEWER = {G.\ A.\ Heuer},
       DOI = {10.1080/0025570x.1985.11977198},
       URL = {https://doi-org.libproxy.clemson.edu/10.1080/0025570x.1985.11977198},
}

@article {etienne1991,
    AUTHOR = {Etienne, Gwihen},
     TITLE = {Tableaux de {Y}oung et solitaire bulgare},
   JOURNAL = {J. Combin. Theory Ser. A},
  FJOURNAL = {Journal of Combinatorial Theory. Series A},
    VOLUME = {58},
      YEAR = {1991},
    NUMBER = {2},
     PAGES = {181--197},
      ISSN = {0097-3165,1096-0899},
   MRCLASS = {05E99 (05E10)},
  MRNUMBER = {1129115},
MRREVIEWER = {Andrea\ Brini},
       DOI = {10.1016/0097-3165(91)90059-P},
       URL = {https://doi-org.libproxy.clemson.edu/10.1016/0097-3165(91)90059-P},
}

@article {griggs1998,
    AUTHOR = {Griggs, Jerrold R. and Ho, Chih-Chang},
     TITLE = {The cycling of partitions and compositions under repeated
              shifts},
   JOURNAL = {Adv. in Appl. Math.},
  FJOURNAL = {Advances in Applied Mathematics},
    VOLUME = {21},
      YEAR = {1998},
    NUMBER = {2},
     PAGES = {205--227},
      ISSN = {0196-8858,1090-2074},
   MRCLASS = {05A17 (11P81)},
  MRNUMBER = {1634705},
MRREVIEWER = {Dennis\ A.\ Eichhorn},
       DOI = {10.1006/aama.1998.0597},
       URL = {https://doi-org.libproxy.clemson.edu/10.1006/aama.1998.0597},
}

@article {popov2005,
    AUTHOR = {Popov, Serguei},
     TITLE = {Random {B}ulgarian solitaire},
   JOURNAL = {Random Structures Algorithms},
  FJOURNAL = {Random Structures \& Algorithms},
    VOLUME = {27},
      YEAR = {2005},
    NUMBER = {3},
     PAGES = {310--330},
      ISSN = {1042-9832,1098-2418},
   MRCLASS = {60J10 (05A99 60C05 91A60)},
  MRNUMBER = {2162601},
MRREVIEWER = {Gh.\ Opri\c san},
       DOI = {10.1002/rsa.20076},
       URL = {https://doi-org.libproxy.clemson.edu/10.1002/rsa.20076},
}

@article {hopkins2006,
    AUTHOR = {Hopkins, Brian and Jones, Michael A.},
     TITLE = {Shift-induced dynamical systems on partitions and
              compositions},
   JOURNAL = {Electron. J. Combin.},
  FJOURNAL = {Electronic Journal of Combinatorics},
    VOLUME = {13},
      YEAR = {2006},
    NUMBER = {1},
     PAGES = {Research Paper 80, 19},
      ISSN = {1077-8926},
   MRCLASS = {05A17 (37B10)},
  MRNUMBER = {2255422},
       DOI = {10.37236/1106},
       URL = {https://doi.org/10.37236/1106},
}

@incollection {hopkins2007garden,
    AUTHOR = {Hopkins, Brian and Sellers, James A.},
     TITLE = {Exact enumeration of {G}arden of {E}den partitions},
 BOOKTITLE = {Combinatorial number theory},
     PAGES = {299--303},
 PUBLISHER = {de Gruyter, Berlin},
      YEAR = {2007},
      ISBN = {978-3-11-019029-8},
   MRCLASS = {05A17 (11P81)},
  MRNUMBER = {2337053},
}

@article {HarrisNguyen2025Depth,
    AUTHOR = {Harris, A. J. and Nguyen, Son},
     TITLE = {Bulgarian solitaire: a new representation for depth generating
              functions},
   JOURNAL = {Australas. J. Combin.},
  FJOURNAL = {The Australasian Journal of Combinatorics},
    VOLUME = {93},
      YEAR = {2025},
     PAGES = {224--273},
      ISSN = {1034-4942,2202-3518},
   MRCLASS = {05A15},
  MRNUMBER = {4982922},
MRREVIEWER = {Warren\ P.\ Johnson},
}

@article {Pham2023Limiting,
    AUTHOR = {Pham, Nhung},
     TITLE = {Limiting behavior in growth of {B}ulgarian solitaire orbits},
   JOURNAL = {Integers},
  FJOURNAL = {Integers. Electronic Journal of Combinatorial Number Theory},
    VOLUME = {23},
      YEAR = {2023},
     PAGES = {Paper No. A58, 25},
      ISSN = {1553-1732},
   MRCLASS = {05A17},
  MRNUMBER = {4628640},
       DOI = {10.5281/zenodo.8214839},
       URL = {https://doi-org.libproxy.clemson.edu/10.5281/zenodo.8214839},
}

@article {akinDavis1985bulgarian,
    AUTHOR = {Akin, Ethan and Davis, Morton},
     TITLE = {Bulgarian solitaire},
   JOURNAL = {Amer. Math. Monthly},
  FJOURNAL = {American Mathematical Monthly},
    VOLUME = {92},
      YEAR = {1985},
    NUMBER = {4},
     PAGES = {237--250},
      ISSN = {0002-9890,1930-0972},
   MRCLASS = {05A17 (90D99)},
  MRNUMBER = {786523},
MRREVIEWER = {A.\ Pultr},
       DOI = {10.2307/2323643},
       URL = {https://doi-org.libproxy.clemson.edu/10.2307/2323643},
}

@article{eriksson1981bulgarisk,
  title={Bulgarisk patiens},
  author={Eriksson, Henrik},
  journal={Elementa},
  volume={64},
  number={4},
  pages={186--188},
  year={1981}
}

@article {canningsHaigh1992montreal,
    AUTHOR = {Cannings, C. and Haigh, J.},
     TITLE = {Montreal solitaire},
   JOURNAL = {J. Combin. Theory Ser. A},
  FJOURNAL = {Journal of Combinatorial Theory. Series A},
    VOLUME = {60},
      YEAR = {1992},
    NUMBER = {1},
     PAGES = {50--66},
      ISSN = {0097-3165,1096-0899},
   MRCLASS = {90D46 (05A99)},
  MRNUMBER = {1156647},
MRREVIEWER = {A.\ Pultr},
       DOI = {10.1016/0097-3165(92)90037-U},
       URL = {https://doi-org.libproxy.clemson.edu/10.1016/0097-3165(92)90037-U},
}

@article {olson2016variants,
    AUTHOR = {Olson, Jeffrey S.},
     TITLE = {Variants of {B}ulgarian solitaire},
   JOURNAL = {Integers},
  FJOURNAL = {Integers. Electronic Journal of Combinatorial Number Theory},
    VOLUME = {16},
      YEAR = {2016},
     PAGES = {Paper No. A8, 11},
      ISSN = {1553-1732},
   MRCLASS = {91A46 (11P81)},
  MRNUMBER = {3458337},
}

@article {erikssonJonssonSjostrand2020limit,
    AUTHOR = {Eriksson, Kimmo and Jonsson, Markus and Sj\"ostrand, Jonas},
     TITLE = {Limit shapes of stable and recurrent configurations of a
              generalized {B}ulgarian solitaire},
   JOURNAL = {Online J. Anal. Comb.},
  FJOURNAL = {Online Journal of Analytic Combinatorics},
    NUMBER = {15},
      YEAR = {2020},
     PAGES = {Paper No. 10, 19},
      ISSN = {1931-3365},
   MRCLASS = {91A99 (05A15)},
  MRNUMBER = {4142376},
MRREVIEWER = {Brian\ Hopkins},
}

@article {MR4011610,
    AUTHOR = {da Silva, Robson and Hopkins, Brian and Sellers, James A.},
     TITLE = {Garden of {E}den states in {A}ustrian solitaire},
   JOURNAL = {European J. Combin.},
  FJOURNAL = {European Journal of Combinatorics},
    VOLUME = {83},
      YEAR = {2020},
     PAGES = {103023, 8},
      ISSN = {0195-6698,1095-9971},
   MRCLASS = {91A46 (05A15 05A19 37B99)},
  MRNUMBER = {4011610},
MRREVIEWER = {Romeo\ Me\v strovi\'c},
       DOI = {10.1016/j.ejc.2019.103023},
       URL = {https://doi.org/10.1016/j.ejc.2019.103023},
}

@article {erikssonJonssonSjostrand2018randomq,
    AUTHOR = {Eriksson, Kimmo and Jonsson, Markus and Sj\"ostrand, Jonas},
     TITLE = {An exponential limit shape of random {$q$}-proportion
              {B}ulgarian solitaire},
   JOURNAL = {Integers},
  FJOURNAL = {Integers. Electronic Journal of Combinatorial Number Theory},
    VOLUME = {18},
      YEAR = {2018},
     PAGES = {Paper No. A58, 22},
      ISSN = {1553-1732},
   MRCLASS = {91A60 (05A17 60J10)},
  MRNUMBER = {3819877},
}

@article {erikssonJonsson2017levels,
    AUTHOR = {Eriksson, Henrik and Jonsson, Markus},
     TITLE = {Level sizes of the {B}ulgarian solitaire game tree},
   JOURNAL = {Fibonacci Quart.},
  FJOURNAL = {The Fibonacci Quarterly. The Official Journal of the Fibonacci
              Association},
    VOLUME = {55},
      YEAR = {2017},
    NUMBER = {3},
     PAGES = {243--251},
      ISSN = {0015-0517,2641-340X},
   MRCLASS = {05A17 (05A15 11P81)},
  MRNUMBER = {3696269},
MRREVIEWER = {Stephan\ G.\ Wagner},
       DOI = {10.1080/00150517.2017.12427770},
       URL = {https://doi-org.libproxy.clemson.edu/10.1080/00150517.2017.12427770},
}

@article{meyerowitzEtAl2026floridian,
  title={Floridian Solitaire: A New Variant of {B}ulgarian Solitaire},
  author={Meyerowitz, Aaron and Curran, Stephen J. and Locke, Stephen C. and Low, Richard M.},
  journal={arXiv:2608.08313},
  year={2026}
}

@incollection {hopkins2009envelopes,
    AUTHOR = {Hopkins, Brian},
     TITLE = {Column-to-row operations on partitions: the envelopes},
 BOOKTITLE = {Combinatorial number theory},
     PAGES = {65--76},
 PUBLISHER = {Walter de Gruyter, Berlin},
      YEAR = {2009},
      ISBN = {978-3-11-020221-2},
   MRCLASS = {05A17 (37E15)},
  MRNUMBER = {2521954},
MRREVIEWER = {Philippe\ Nadeau},
}

@article{bojanov1981problem,
  title={Problem solution 4},
  author={Bojanov, Borislav},
  journal={Obuchenieto po matematika},
  volume={24},
  number={5},
  pages={59--60},
  year={1981}
}

@article {yeh1995remarkable,
    AUTHOR = {Yeh, Yeong Nan},
     TITLE = {A remarkable endofunction involving compositions},
   JOURNAL = {Stud. Appl. Math.},
  FJOURNAL = {Studies in Applied Mathematics},
    VOLUME = {95},
      YEAR = {1995},
    NUMBER = {4},
     PAGES = {419--432},
      ISSN = {0022-2526,1467-9590},
   MRCLASS = {05A17 (90D42)},
  MRNUMBER = {1358813},
MRREVIEWER = {Richard\ J.\ Nowakowski},
       DOI = {10.1002/sapm1995954419},
       URL = {https://doi.org/10.1002/sapm1995954419},
}

\end{document}